%% file: agt-2-20.tex
\documentclass{gtart}

\input agtout

\lognumber{20}
\volumenumber{2}
\volumeyear{2002}
\papernumber{20}
\published{25 May 2002}
\pagenumbers{403}{432}
\received{15 March 2002}
\revised{9 May 2002}
\accepted{15 May 2002}

\usepackage{epsf}
\usepackage{amssymb}
\usepackage{amsmath}

\usepackage{xy}
\xyoption{all}

\def\C{{\Bbb C}}
\def\Z{{\Bbb Z}}

\def\CP{{\Bbb{CP}}}

\def\P{{\Bbb P}}
\def \a{\alpha}
\def \be{\beta}
\def \Dl{\Delta}

\def \g{\gamma}
\def \G{\Gamma}

\def \s{\sigma}

\def \la{\langle}

\def \sbt{\subset}

\def \ra{\rangle}
\def \1{^{-1}}
\def \2{^{-2}}

\def\Pitil{\widetilde\Pi}

\def \Aff{\operatorname{Aff}}

\def \Gal{\operatorname{Gal}}

\def\GG{\Pitil_1}

\newcommand\semidirect{\ltimes}
\newcommand\Ker{\operatorname{Ker}}
\renewcommand\Im{\operatorname{Im}}
\newcommand\sg[1]{{\left<{#1}\right>}}
\newcommand\FIGURE[3][]{{\begin{figure}[ht!]\epsfysize=0cm \cl{\epsfbox{\figs{#2}}}\caption{#1}\label{#3}\end{figure}}}
\newcommand\FIGUREx[4][]{{\begin{figure}[ht!]\epsfxsize=#4 \cl{\epsfbox{\figs{#2}}}\caption{#1}\label{#3}\end{figure}}}
\newcommand\FIGtab[1]{{$\stackrel{}{\hbox{\epsfbox{\figs{#1}}}}$}}
\newcommand\FIGtabx[2]{{$\stackrel{}{\hbox{\epsfxsize=#2\epsfbox{\figs{#1}}}}$}}
\def\st{{such that}}
\newcommand\set[1]{{\{#1\}}}
\newcommand\eq[1]{{(\ref{#1})}}
\newcommand\eqs[2]{{(\ref{#1})--(\ref{#2})}}
\newcommand\Trip[2]{{{#1}{#2}{#1}={#2}{#1}{#2}}}
\def\limiits{{}} 
\newcommand\pat[1]{{\bf #1}}
\newcommand\tilpat[1]{{\bf\tilde{#1}}}

\newcommand\figs[1]{#1}

\newtheorem{theorem}{Theorem}[section] 

\newtheorem{tab}[theorem]{Table}

\newtheorem{corollary}[theorem]{Corollary}
\newtheorem{lemma}[theorem]{Lemma}

\newtheorem{defn}[theorem]{Definition}
\newtheorem{proposition}[theorem]{Proposition}
\newtheorem{remark}[theorem]{Remark}

\newenvironment{Relax}{\relax}{\relax}
\let\Bbb\mathbb

\begin{document}

\title{The fundamental group of a Galois cover of $\C\P^1 \times T$}
\asciititle{The fundamental group of a Galois cover of CP^1 X T}
\covertitle{The fundamental group of a Galois cover of ${\noexpand
\bf CP}^1 \times T$}

\author{Meirav Amram\\David Goldberg\\Mina Teicher\\Uzi Vishne}

\shortauthors{Amram, Goldberg, Teicher and Vishne}

\address{Meirav Amram and Mina Teicher:\\Department of Mathematics, 
Bar-Ilan University, Ramat-Gan 52900, Israel\\\smallskip\\David Goldberg:\\Department of 
Mathematics, Colorado State University\\Fort Collins, CO 80523-1874, USA\\\smallskip\\Uzi Vishne:\\Einstein Institute of Mathematics, Givat Ram Campus\\The 
Hebrew University of Jerusalem, Jerusalem 91904, Israel}

\asciiaddress{Meirav Amram and Mina Teicher:\\Department of Mathematics, 
Bar-Ilan University, Ramat-Gan 52900, Israel\\David Goldberg:\\Department of 
Mathematics, Colorado State University\\Fort Collins, CO 80523-1874, USA\\Uzi Vishne:\\Einstein Institute of Mathematics, Givat Ram Campus\\The 
Hebrew University of Jerusalem, Jerusalem 91904, Israel}

\email{amram@mi.uni-erlangen.de, meirav@macs.biu.ac.il,
goldberg@math.colostate.edu,
teicher@macs.biu.ac.il,
vishne@math.huji.ac.il}

\begin{abstract}
Let $T$ be the complex projective torus, and $X$ the surface
$\C\P^1 \times T$. Let $X_{\Gal}$ be its Galois cover with respect
to a generic projection to $\C\P^2$. In this paper we compute the
fundamental group of $X_{\Gal}$, using the degeneration and
regeneration techniques, the Moishezon-Teicher braid monodromy
algorithm and group calculations. We show that $\pi_1(X_{\Gal}) =
\Z^{10}$.
\end{abstract}

\asciiabstract{Let T be the complex projective torus, and X the
surface CP^1 X T. Let X_Gal be its Galois cover with respect to a
generic projection to CP^2. In this paper we compute the fundamental
group of X_Gal, using the degeneration and regeneration
techniques, the Moishezon-Teicher braid monodromy algorithm and group
calculations. We show that pi_1(X_Gal) = Z^10.}



\maketitle

\begin{Relax}\end{Relax}

\section{Overview}\label{sec:1}
Let $T$ be a complex torus in $\C\P^2$. We compute the fundamental
group of the Galois cover with respect to a generic map of the
surface $X = \C \P^1 \times T$ to $\C\P^2$. We embed $X$ into a
projective space using the Segre map $\C\P^1\times\C\P^2\to\C\P^5$
defined by $(s_0,s_1)\times(t_0,t_1,t_2)\mapsto
(s_0t_0,s_1t_0,s_0t_1,s_1t_1,s_0t_2,s_1t_2)$. Then, a generic
projection $f \co X \rightarrow \C\P^2$ is obtained by projecting
$X$ from a general plane in $\C\P^5-X$ to $\C\P^2$.
The Galois cover can now be defined as the closure of the $n$-fold
fibered product
$X_{\Gal}=\overline{X\times_f\cdots\times_fX-\Delta}$ where $n$ is
the degree of the map $f$, and $\Delta$ is the generalized
diagonal.  The closure is necessary because the branched fibers
are excluded when $\Delta$ is omitted.

The fundamental group $\pi_1(X_{\Gal})$ is related to the fundamental group 
of the complement of the branch curve. The latter is an important invariant
of $X$, which can be used to classify algebraic surfaces of a general
type, up to deformations. Such an invariant is finer 
than the famous Sieberg-Witten invariants 
and thus can serve as a tool to distinguish diffeomorphic 
surfaces which are not deformation of each other 
(see \cite{KuTe}), \cite{DD1} and \cite{DD2})
a problem which is referred to as the 
Diff-Def problem.
The algorithms and problems that arise in the computation of these two 
types of groups are related, and one hopes to be able to compute such 
groups for various types of surfaces.

Since the induced map $X_{\Gal}\to\C\P^2$ has the same branch
curve $S$ as $f\colon X\to\C\P^2$, the fundamental group
$\pi_1(X_{\Gal})$ is related to $\pi_1(\C\P^2-S)$. In fact it is a
normal subgroup of the quotient of $\pi_1(\C\P^2-S)$ by the normal
subgroup generated by the squares of the standard generators. In
this paper we employ braid monodromy techniques, the van Kampen theorem
and various computational methods of groups to compute a
presentation for the quotient $\Pitil_1$ from which
$\pi_1(X_{\Gal})$ can be derived. Our main result is that
$\pi_1(X_{\Gal}) = \Z^{10}$ (Theorem \ref{main}).

The fundamental group of the Galois cover of the surface 
$X = \C \P^1 \times T$ is a step in computing the same group for 
$T \times T$ \cite{TxT}, and will later appear in local 
computations of fundamental groups of the Galois cover of $K3$-surfaces.

It turns out that 
a property of the Galois covers that 
were treated before (see \cite{MoRoTe}, \cite{MoTe1} or \cite{MoTe3})
is lacking in the Galois cover of $X = \C\P^1\times T$.
In all the cases computed so far, $X$ had the
property that the fundamental group of the graph defined on the
planes of the degenerated surface $X_0$ by connecting every two
intersecting planes, is generated by the cycles around the
intersection points. Our surface, together with a parallel work on
$T\times T$ \cite{TxT}, are the first cases for which this
assumption does not hold. The significance of this 'redundancy'
property of $S_0$ will be explained in Section \ref{sec:6} (and in
more details in \cite{TxT}).

The paper is organized as follows. In Section \ref{sec:2} we
describe the degeneration of the surface $X = \C\P^1\times T$ and
the degenerated branch curve. In Sections \ref{sec:3} and
\ref{sec:4} we study the regeneration of this curve and its braid
monodromy factorization. We also get a presentation for
$\pi_1(\C^2-S,u_0)$, the fundamental group of the complement of
the regenerated branch curve in $\C^2$, see Theorem \ref{pres}. In
Section \ref{sec:5} we present the homomorphism
$\psi\colon\tilde{\pi}_1(\C^2-S,u_0)\to S_6$, whose kernel is
$\pi_1(X_{\Gal}^{\Aff})$. In Section \ref{sec:6} we study a
natural Coxeter quotient of $\pi_1(\C^2-S,u_0)$, and give its
structure. In Sections \ref{sec:7} and \ref{sec:8} we use the
Reidmeister-Schreir method to give a new presentation for
$\pi_1(X_{\Gal}^{\Aff})$, see Theorem \ref{th100}. In Section
\ref{sec:9} we introduce the projective relation, and prove the
main result about the structure of  $\pi_1(X_{\Gal})$.

\rk{Acknowledgements} Meirav Amram is partially supported by the Emmy
Noether Research Institute for Mathematics (center of the Minerva
Foundation of Germany), the Excellency Center ``Group Theoretic
Methods in the Study of Algebraic Varieties" of the Israel Science
Foundation, and EAGER (EU network, HPRN-CT-2009-00099).

David Goldberg is partially supported by the Emmy Noether Research
Institute for Mathematics, and the Minerva Foundation of Germany. 

Uzi Vishne is partially supported by the Edmund Landau Center for
Research in Mathematical Analysis and Related Subjects.

\section{Degeneration of  $\C\P^1\times T$}\label{sec:2}
In the computation of braid monodromies it is often useful to
replace the surface with a degenerated object, made of copies of
$\CP^2$. It is easy to see that $T$ degenerates to a triangle of
complex projective lines (see \cite[Subsection 1.6.3]{Am}), so $X$
degenerates to a union of three quadrics $Q_1$, $Q_2$, and $Q_3$,
$Q_i = \CP^1\times \CP^1$, which we denote by $X_1$, see Figure
\ref{threeSQ}.

\FIGURE[{The space $X_1$}]{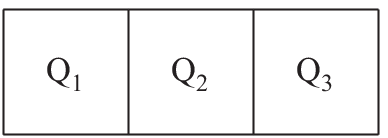}{threeSQ} 

Each square in Figure \ref{threeSQ} represents a quadric surface.
Since $T$ degenerates to a triangle, $Q_1$ and $Q_3$ intersect, so
the left and right edges of $X_1$ are identified. Therefore, we
can view $X_1$ as a triangular prism.

Each quadric in $X_1$ can be further degenerated to a union of two
planes.  In Figure \ref{sixtri} this is represented by a diagonal
line which divides each square into two triangles, each one
isomorphic to $\CP^2$.

\FIGURE[The simplicial complex $X_0$]{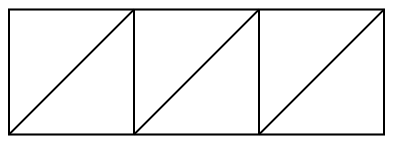}{sixtri}

We shall refer to this diagram as the simplicial complex of $X_0$.
A common edge between two triangles represents the intersection
line of the two corresponding planes.  The union of the
intersection lines is the ramification locus in $X_0$ of
$f_0\colon X_0\to \C\P^2$, denoted by $R_0$. Let  $S_0 = f_0(R_0)$
be the degenerated branch curve. It is a line arrangement,
composed of all the intersection lines.

A vertex in the simplicial complex represents an intersection
point of three planes.  The vertices represent singular points of
$R_0$. Each of these vertices is called a 3-point (reflecting the
number of planes which meet there).

The vertices may be given any convenient enumeration.  We have chosen
left to right, bottom to top enumeration, see Figure \ref{sixp}.
The extreme vertices are pairwise identified,
as well as the left and right edges.

{\nocolon\FIGUREx{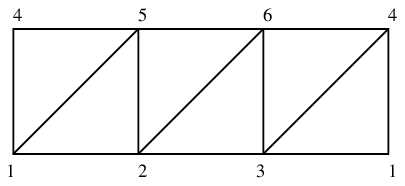}{sixp}{5.2cm}}

We create an enumeration of the edges based upon the enumeration
of the vertices using reverse lexicographic ordering: if $L_1$ and
$L_2$ are two lines with end points $\alpha _1, \beta _1$ and
$\alpha _2, \beta _2$ respectively $(\alpha _1 < \beta _1, \alpha
_2 < \beta _2)$, then $L_1 < L_2$ iff $\beta_1 < \beta _2$, or
$\beta _1 = \beta _2$ and $\alpha _1 < \alpha _2$. The resulting
enumeration is shown in Figure \ref{sixe}. This enumeration
dictates the order of the regeneration of the lines to curves, see
the next section. The horizontal lines at the top and bottom do
not represent intersections of planes and hence are not numbered.

{\nocolon\FIGUREx{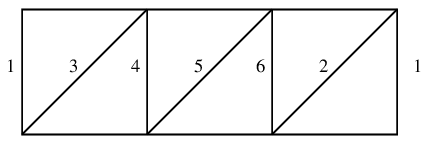}{sixe}{5.5cm}}

We enumerate the triangles $\set{P_i}_{i=1}^6$ also according to
the enumeration of vertices in reverse lexicographic order. If
$P_i$ and $P_j$ have vertices $\alpha_1, \alpha_2, \alpha_3$ and
$\beta_1, \beta_2, \beta_3$ respectively, with $\alpha_1 <
\alpha_2 < \alpha_3$ and $\beta_1 < \beta_2 < \beta_3$, then $P_i
< P_j$ iff $\alpha_3 < \beta_3$, or $\alpha_3 = \beta_3$ and
$\alpha_2 < \beta_2$, or $\alpha_3 = \beta_3$, $\alpha_2 =
\beta_2$ and $\alpha_1 < \beta_1$. The enumeration is shown in
Figure~\ref{sixinside}.

{\nocolon\FIGUREx{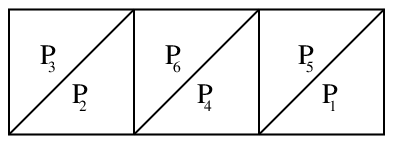}{sixinside}{5cm}}

\section{Regeneration of the Branch Curve}\label{sec:3}

\subsection{The Braid Monodromy of $S_0$}\label{ssec:31}
Starting  from the branch curve $S_0$, we reverse the steps in the
degeneration of $X$ to regenerate the braid monodromy of $S$.
Figure \ref{dege} shows the three steps to recover the original
object $X_3=X$.

\FIGUREx[The regeneration process]{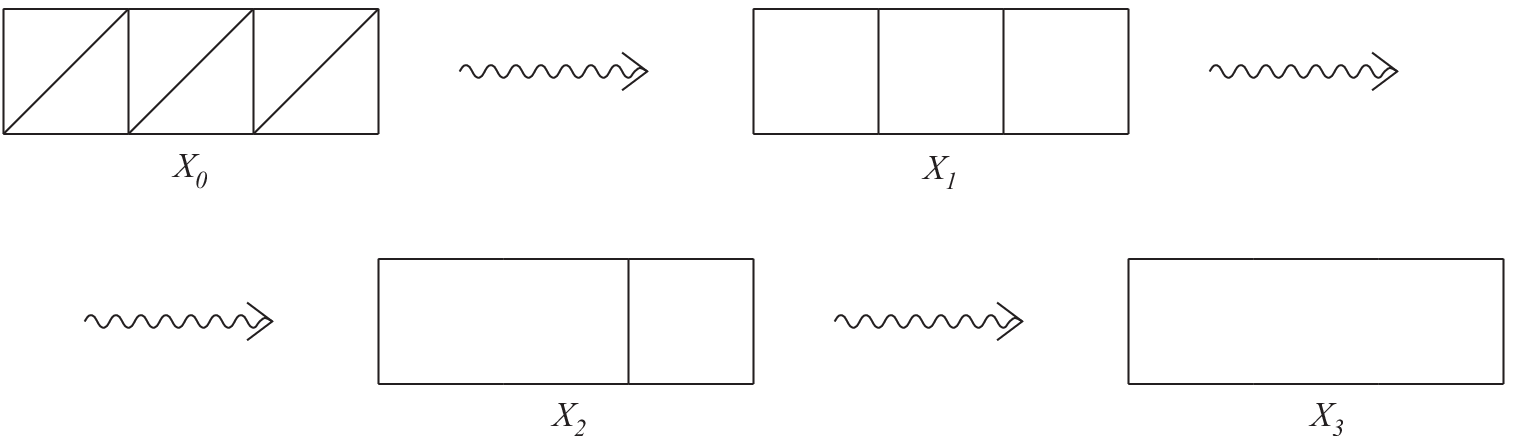}{dege}{11cm}

Recall that $X$ comes with an embedding to $\CP^5$.    At each
step of the regeneration, the generic projection $\C\P^5\to\C\P^2$
restricts to a generic map $f_i\colon X_i\to\C\P^2$. Let $R_i\sbt
X_i$ be the ramification locus of $f_i$ and $S_i\sbt\C\P^2$ the
corresponding branch curve.

We have enumerated the six planes $P_1,\ldots,P_6$ which comprise
$X_0$, the six intersection lines $\hat L_1,\ldots,\hat L_6$, and
their six intersection points $\hat V_1,\ldots,\hat V_6$.  Let
$L_i$ and $V_j$ denote the projections of $\hat L_i$ and $\hat
V_j$ to $\C\P^2$ by the map $f_0$.  Clearly
$R_0=\bigcup\limiits_{i=1}^6\hat L_i$ 
and $S_0=\bigcup \limiits_{i=1}^6L_i$.  Let $C$ be the line arrangement
consisting of all lines through pairs of the $V_j$s.  The
degenerated branch curve $S_0$ is a sub-arrangement of $C$.  Since
$C$ is a dual to generic arrangement, Moishezon's results in \cite{Mo} 
(and later on Theorem IX.2.1 in
\cite{MoTe2}) gives us a braid monodromy factorization for $C$:
$\Dl_C^2=\prod \limiits_{j=1}^6 C_j \Dl_j^2$ where $\Dl_j^2$ is the
monodromy around $V_j$ and the $C_j$ consist of products of the
monodromies around the other intersections points of $C$.  This
factorization can be restricted to $S_0$ by removing from the
braids all strands which correspond to lines of $C$ that do not
appear in $S_0$, and deleting all factors which correspond to
intersections in $C$ that do not appear in $S_0$.  Thus we get a
braid monodromy factorization: $\Dl_{S_0}^2=\prod \limiits_{j=1}^6
\tilde{C}_j \tilde{\Dl}_j^2$.  The $\tilde{C}_j$ and
$\tilde{\Dl}_j^2$ and their regenerations are formulated more
precisely in the following subsections.

\subsection{$\tilde{\Delta}^2_j$ and its Regeneration}\label{ssec:33} 

Consider an affine piece of $S_0\sbt\C\P^2$ and take a generic
projection $\pi\colon\C^2\to\C$.  Let $N$ be the set of the projections of the singularities and branch points with respect to $\pi$. 
Choose $u\in\C-N$ and let
$\C_u=\pi^{-1}(u)$ be a generic fiber.

A path from a point $j$ to a point $k$ below 
the real line is denoted by $\underline{\pat{z}}_{jk}$, and the corresponding
halftwist by $\underline{\pat{Z}}_{jk}$.

Two lines $\hat L_j$ and $\hat L_k$ which meet in $X_0$ give rise to
braids connecting $j=L_j\cap\C_u$ and $k=L_k\cap\C_u$, namely a fulltwist
$\underline{\pat{Z}}_{jk}^2$
of $j$ and $k$. This is done in the following way: 
let $V_i = L_j \cap L_k$ be the intersection point, then
$\tilde\Dl_i^2=\underline{\pat{Z}}_{jk}^2$.

We shall analyze the regeneration of the local braid monodromy of $S$
in a small neighborhood of each $V_i$. The case of non-intersecting lines 
(which give 'parasitic intersections') is discussed in the next subsection.

\begin{figure}[ht!]
\begin{center}
\begin{tabular}{cccccc} 
	{}\FIGtab{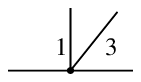}
	& \FIGtab{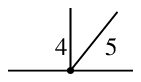}
	& \FIGtab{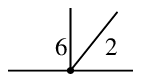}
	& \FIGtab{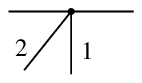}
	& \FIGtab{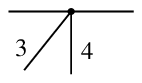}
	& \FIGtab{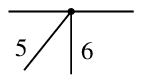}
\end{tabular}
\end{center}
\caption{Enumeration around the 3-points}\label{3pt}
\end{figure}

The degenerated branch curve $S_0$ has six singularities coming 
from the 3-points of $R_0$, shown in Figure~\ref{3pt}.

Each pair of lines
intersecting at a 3-point regenerates in $S_1$ (the branch curve of  $X_1$) 
as follows: the
diagonal line becomes a conic, and the 
vertical line is tangent to it. In
the next step of the regeneration the point of tangency becomes
three cusps according to the third regeneration rule (which was quoted in 
\cite{Mo} and proven in 
\cite[p.337]{MoTe4}).  This is enough
information to compute $H_{V_i}$, the local braid monodromy
of $S$ in a neighborhood of $V_i$, see these specific computations for this case in  \cite[Subsection 1.10.4]{Am}. 

In the regeneration, each point $\alpha$ on the typical fiber is 
replaced by two close points $\alpha, \alpha'$. 
Denote by $\rho_{\a} = \pat{Z}_{\a\a'}$ the counterclockwise 
halftwist of $\a$ and $\a'$.

The following table presents the global form of the local braid monodromies, as quoted in \cite{Mo}, and presents also the application of this global form to our case.

\begin{figure}[ht!]
\cl{\epsfbox{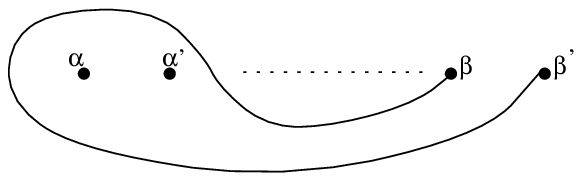}\qquad\quad\epsfbox{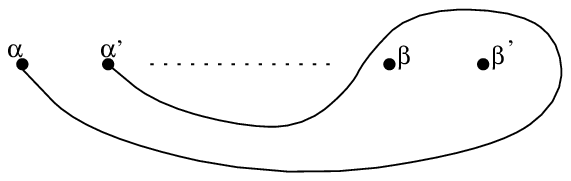}}
\caption{$\pat{Z}_{\be\be'(\a)}$ for $i = 1,2,4$\qquad\qquad\quad Figure 9:
$\pat{Z}_{\a\a'(\be)}$ for $i = 3,5,6$}\label{fig10}
\end{figure}

\begin{tab}\label{HV}
The local braid monodromies $H_{V_i}$ are as follows. For every fixed $i$,
let $\a < \be$ be the lines intersecting at $V_i$. 
Let $Z^{3}_{\a\a',\be} = (\underline{\pat{Z}}_{\a'\be}^3)_{\rho_{\a}}\cdot 
\underline{\pat{Z}}_{\a'\be}^{3} \cdot
(\underline{\pat{Z}}^{3}_{\a'\be})_{\rho_{\a}^{-1}}$ 
and 
$Z^{3}_{\a',\be\be'} = (\underline{\pat{Z}}_{\a'\be}^3)_{\rho_{\be}}\cdot 
\underline{\pat{Z}}_{\a'\be}^{3}
\cdot (\underline{\pat{Z}}^{3}_{\a'\be})_{\rho_{\be}^{-1}}$.

For $i=1,2,4$ we have 
$H_{V_i} = \pat{Z}^{3}_{\a\a',\be'}\cdot\pat{Z}_{\be\be'(\a)}$, where
$\pat{Z}_{\be\be'(\a)}$ is the halftwist corresponding to the path
shown in Figure~\ref{fig10}.
For $i=3,5,6$ we have 
$H_{V_i} = \pat{Z}^{3}_{\a',\be\be'}\cdot\pat{Z}_{\a\a'(\be)}$, where
$\tilpat{Z}_{\a\a'(\be)}$ is the halftwist corresponding to the path 
shown in Figure~9.

In paricular we have
\begin{eqnarray*}
H_{V_1} & = & \pat{Z}_{11',3'}^{3}\cdot \pat{Z}_{33'(1)}, \\
H_{V_2} & = & \pat{Z}_{44',5'}^{3}\cdot \pat{Z}_{55'(4)}, \\
H_{V_3} & = & \pat{Z}_{2',66'}^{3}\cdot \pat{Z}_{22'(6)},\\
H_{V_4} & = & \pat{Z}_{11',2'}^{3}\cdot \pat{Z}_{22'(1)}, \\
H_{V_5} & = & \pat{Z}_{3',44'}^{3}\cdot \pat{Z}_{33'(4)},\\
H_{V_6} & = & \pat{Z}_{5',66'}^{3}\cdot \pat{Z}_{55'(6)}.
\end{eqnarray*}
\end{tab}

\addtocounter{figure}{1}
\def\strut{\vrule width 0pt height 12pt depth 5pt}

\begin{figure}\small
\begin{center}
\begin{tabular}{|l|c|l|l|} \hline
    pt.
    & \strut the braid
    & exp.
    & the path representing the braid
\\
\hline
    {}
    & \strut$\rho_1^i  \underline{\pat{Z}}_{1'3'} \rho_1^{-i}\quad i=0,1,2$
    & 3
    &
\\
[-.5cm]
    $V_1$
    & \strut$\pat{Z}_{33'(1)}$
    & 1
    & \FIGtabx{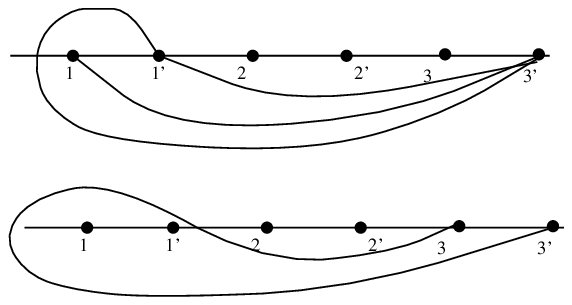}{1.8in}
\\
\hline
    {}
    & \strut$\rho_4^i\underline{\pat{Z}}_{4'5'} \rho_4^{-i}\quad i=0,1,2$
    & 3
    &
\\
[-.5cm]
    $V_2$
    & \strut$\pat{Z}_{55'(4)}$
    & 1
    & \FIGtabx{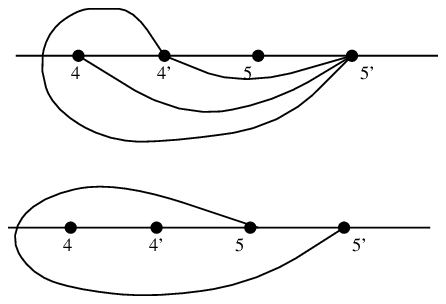}{1.5in}
\\
\hline
    {}
    & \strut$\rho_6^i\underline{\pat{Z}}_{2'6'} \rho_6^{-i}\quad i=0,1,2$
    & 3
    &
\\
[-.5cm]
    $V_3$
    & \strut$\pat{Z}_{22'(6)}$
    & 1
    & \FIGtabx{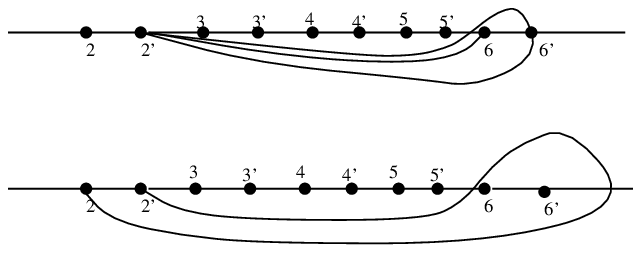}{2.3in}
\\
\hline
    {}
    & \strut$\rho_1^i\underline{\pat{Z}}_{1'2'} \rho_1^{-i}\quad i=0,1,2$
    & 3
    &
\\
[-.5cm]
    $V_4$
    & \strut$\pat{Z}_{22'(1)}$
    & 1
    & \FIGtabx{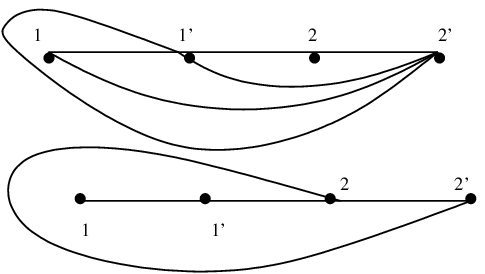}{1.8in}
\\
\hline
    {}
    & \strut$\rho_4^i\underline{\pat{Z}}_{3'4} \rho_4^{-i}\quad i=0,1,2$
    & 3
    &
\\
[-.5cm]
    $V_5$
    & \strut$\pat{Z}_{33'(4)}$
    & 1
    & \FIGtabx{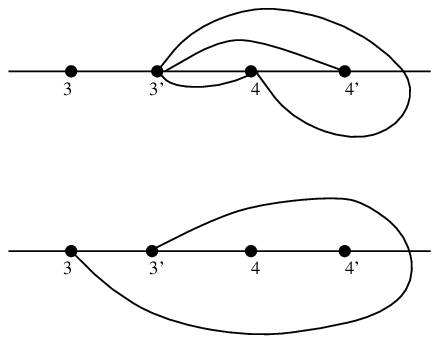}{1.5in}
\\
\hline
    {}
    & \strut$\rho_6^i\underline{\pat{Z}}_{5'6} \rho_6^{-i}\quad i=0,1,2$
    & 3
    &
\\
[-.5cm]
    $V_6$
    & \strut$\pat{Z}_{55'(6)}$
    & 1
    & \FIGtabx{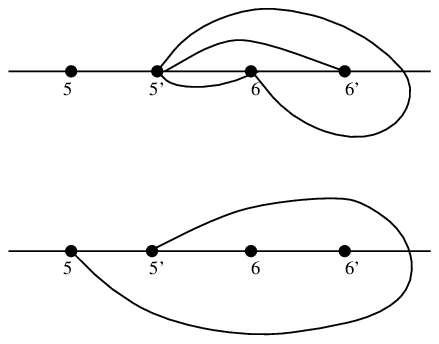}{1.5in}
\\
\hline
\end{tabular}
\end{center}
\caption{Monodromy factorizations}\label{Vtab}
\end{figure}

The table given in Figure \ref{Vtab} presents the six monodromy
factorization, one for every point $V_1,\dots,V_6$. For each
point, the first path represents three factors obtained from
cusps, and the other represents the fourth factor, obtained from
the branch point (as shown in Figures \ref{fig10} and 9). 
The relations obtained from these braids are given
in Theorem \ref{pres}.

\subsection{$\tilde{C}_j$ and its Regeneration}\label{ssec:32} 

There are lines which do not meet in $X_0$ but whose images
meet in $\C^2$.  Such an intersection is called a parasitic
intersection. Each pair of disjoint lines $\hat L_i$ and $\hat L_j$ 
give rise to a certain fulltwist, see \cite[Theorem IX.2.1]{MoTe2}.
This is denoted as $\tilpat{Z}_{ij}^2$, corresponding to a path 
$\tilde{z}_{ij}$, running from $i$ over the points up to $j_0$, then under 
$j_0$ up to $j$, where $j_0$ is the least numbered line which shares the 
upper vertex of $L_j$.

\begin{figure}[ht!]\small
\begin{center}
\begin{tabular}{|c|c|c|c|} \hline
factor & corresponding paths & singularity types & degrees \\ [.2cm]\hline
    $\tilpat{Z}_{ii',jj'}$
    &   \FIGtab{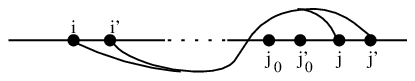}
    & four nodes & 2,2,2,2 \\ [.2cm]\hline
\end{tabular}
\end{center}
\nocolon\caption{}\label{4parts}
\end{figure}

As discussed in \cite{MoTe1} and \cite{Am} the degree of the
regenerated branch curve $S$ is twice the degree of $S_0$.
Consequently each line $L_i\sbt S_0$ divides locally into two
branches of $S$ and each $i=L_i\cap\C_u$ divides into two points,
$i$ and $i'$. According to the second regeneration rule
(quoted in \cite{Mo} and proved in \cite[p.~337]{MoTe4}) the fulltwist  
$\tilpat{Z}_{ij}^2$ 
becomes  $\tilpat{Z}_{ii',jj'}^2$, which compounds four nodes of $S$, 
namely $\tilpat{Z}_{ij}^2$, $\tilpat{Z}_{ij'}^2$,
$\tilpat{Z}_{i'j}^2$ and $\tilpat{Z}_{i'j'}^2$ as shown in 
Figure~\ref{4parts}. 
These are the factors in the
regenerations $C'_j$ of the $\tilde C_j$.
In Table \ref{Ctab} we construct the paths which 
correspond to these braids in our case. 

Pick a base point $u_0\in\C_u$ in the generic fiber
$\C_u=\pi^{-1}(u)$.  The fundamental group $\pi_1(\C_u-S,u_0)$ is
freely generated by $\{\G_j,\G_{j'}\}_{j=1}^6$, where $\G_j$ and
$\G_{j'}$ are loops in $\C_u$ around $j$ and $j'$ respectively.
Let us explain how to create such generators. Define the
generators $\G_j$ and $\G_j'$ in two steps.  First select a path
$\alpha_j$ from $u_0$ to a point $u_j$ close to $j$ and $j'$. Next
choose small counterclockwise circles $\eta_j$ and $\eta'_j$ starting from $u_j$ and
circling $j$ and $j'$ respectively.  Use these three paths to
build the generators $\G_j=\alpha_j\eta_j\alpha_j^{-1}$ and
$\G_{i'}=\alpha_j\eta'_j\alpha_j^{-1}$.  

By the van Kampen Lemma \cite{vK}, there is a surjection from
$\pi_1(\C_u-S,u_0)$ onto  $\Pi_1=\pi_1(\C^2-S,u_0)$. The images of
$\{\G_j,\G_{j'}\}_{j=1}^6$ generate $\Pi_1$. By abuse of notation
we denote them also by  $\{\G_j,\G_{j'}\}_{j=1}^6$. By the van
Kampen Theorem \cite{vK}, each braid in the braid monodromy
factorization of $S$ induces a relation on $\Pi_1$ through its
natural action on $\C_u-\{j,j'\}_{j=1}^6$  \cite{Mo}. A
presentation for
\begin{equation}\label{Pitildef}
\Pitil_1=\frac{\pi_1(\C^2-S,u_0)}{\sg{\G_j^2,\G_{j'}^2}}
\end{equation}
is thus immediately obtained from a presentation of $\Pi_1$. 

The algorithm 
used to compute a relation from a braid
is explained in \cite[Section 0.7]{MoTe1}, see also
\cite[Section 1.11]{Am}. 

Moishezon claimed in \cite{Mo} that the braid
monodromy factorization is invariant under complex conjugation of
$\C_u$. Later it was proven in Lemma 19 of \cite{MoTe4}.
Therefore we can include the complex conjugate paths and relations in
the table.
 For simplicity of notation we will use the following shorthand:
$\G_{ii'}$ will stand for either $\G_i$ or $\G_{i'}$; $\G_{\check
ii'}$ will stand for either $\G_{i'}\G_i\G_{i'}$ or $\G_{i'}$;
$\G_{j\hat j}$ will stand for either $\G_j$ or $\G_j\G_{j'}\G_j$.
\smallskip

\begin{tab}\label{Ctab} 
We present the relations induced by 
the van Kampen Theorem from the paths, one for every pair of non-intersecting 
lines $\hat L_i, \hat L_j$.
{\ }
\begin{itemize}
\item $\tilpat{Z}_{22',33'}^2$.
$\vcenter{\hbox{\epsfxsize=9.75cm \epsfbox{\figs{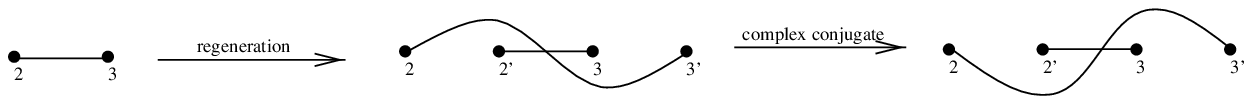}}}}$

The relations: $[\G_{\check{2}2'}, \G_{33'}] = 1$ (from the path
itself), and $[\G_{22'},\G_{3 \hat{3}}] = 1$ (from the complex
conjugate).

\item $\tilpat{Z}_{11',44'}^2$.
$\vcenter{\hbox{\epsfxsize=9.75cm \epsfbox{\figs{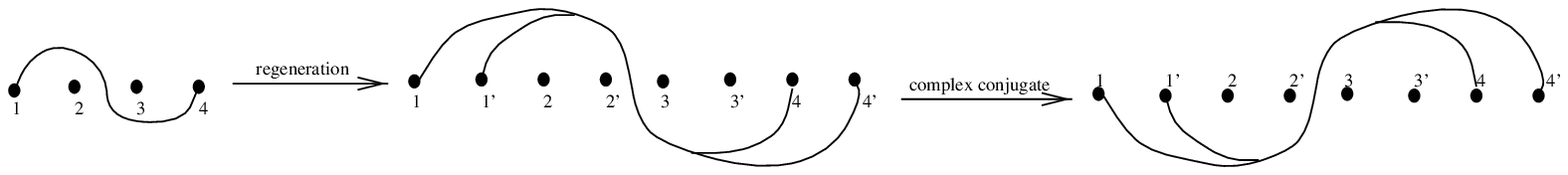}}}}$

Relations: $[\G_{2'}\G_2\G_{\check{1}1'} \G_2 \G_{2'}, \G_{44'}] =
1$ and $[\G_{11'} , \G_3 \G_{3'} \G_{4 \hat{4}} \G_{3'}\G_{3}]
=1$.

\item $\tilpat{Z}_{22',44'}^2$.
$\vcenter{\hbox{\epsfxsize=9.75cm \epsfbox{\figs{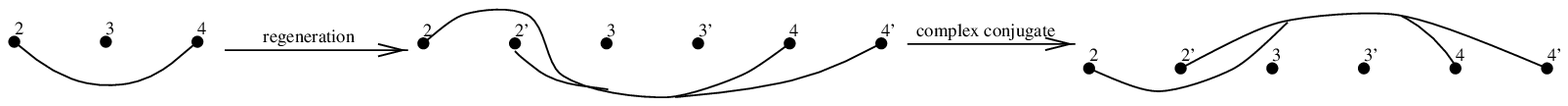}}}}$

$[\G_{\check{2}2'}, \G_{44'}] = 1$ and $[\G_{22'} , \G_3 \G_{3'}
\G_{4 \hat{4}} \G_{3'}\G_{3}] =1$. \vspace{-1cm}

\item  $\tilpat{Z}_{11',55'}^2$.
$\vcenter{\vspace{1cm}\hbox{\epsfxsize=9.75cm\epsfbox{\figs{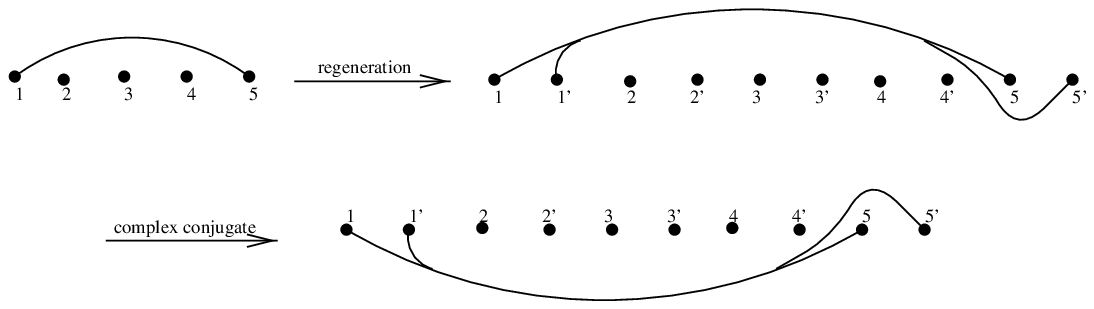}}}}$

$[\G_{4'}\G_{4}\G_{3'}\G_{3}\G_{2'}\G_{2}\G_{\check{1}1'}\G_{2}\G_{2'}\G_{3}\G_{3'}\G_{4}\G_{4'},
\G_{55'}] = 1$ and $[\G_{11'} , \G_{5 \hat{5}} ] =1$.
\vspace{-1cm}

\item $\tilpat{Z}_{22',55'}^2$.
$\vcenter{\vspace{1cm}\hbox{\epsfbox{\figs{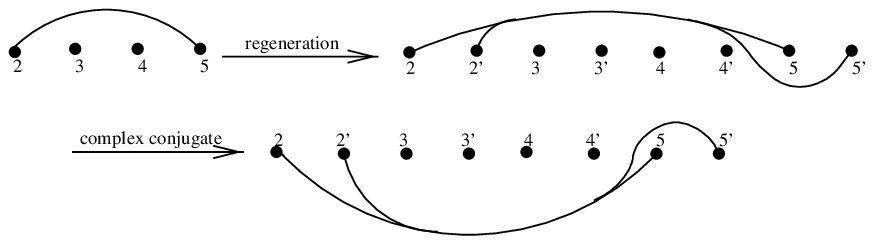}}}}$

$[\G_{4'}\G_{4}\G_{3'}\G_{3}\G_{\check{2}2'}\G_{3}\G_{3'}\G_{4}\G_{4'},
\G_{55'}] = 1$ and $[\G_{22'} , \G_{5 \hat{5}} ] =1$.
\vspace{-0.7cm}

\item $\tilpat{Z}_{33',55'}^2$.
$\vcenter{\vspace{1cm}\hbox{\epsfbox{\figs{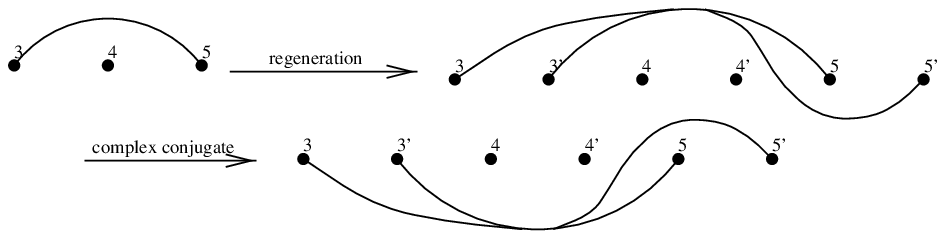}}}}$

$[\G_{4'}\G_{4}\G_{\check{3}3'} \G_4 \G_{4'}, \G_{55'}] = 1$ and
$[\G_{33'} , \G_{5 \hat{5}} ] =1$. \vspace{-1cm}

\item $\tilpat{Z}_{11',66'}^2$.
$\vcenter{\vspace{1cm}\hbox{\epsfxsize=9.75cm\epsfbox{\figs{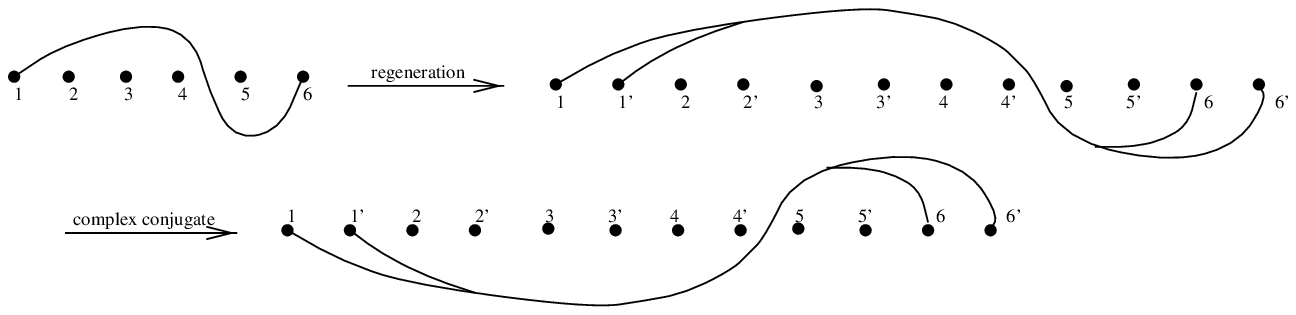}}}}$

Relations: 
$[\G_{4'}\G_{4}\G_{3'}\G_{3}\G_{2'}\G_{2}\G_{\check{1}1'}\G_{2}\G_{2'}\G_{3}\G_{3'}\G_{4}\G_{4'},\G_{66'}] = 1$\nl
and 
$[\G_{11'},\G_5 \G_{5'} \G_{6 \hat{6}} \G_{5'}\G_5]=1$. \vspace{-1cm}

\item $\tilpat{Z}_{33',66'}^2$.
$\vcenter{\vspace{1cm}\hbox{\epsfxsize=9.75cm\epsfbox{\figs{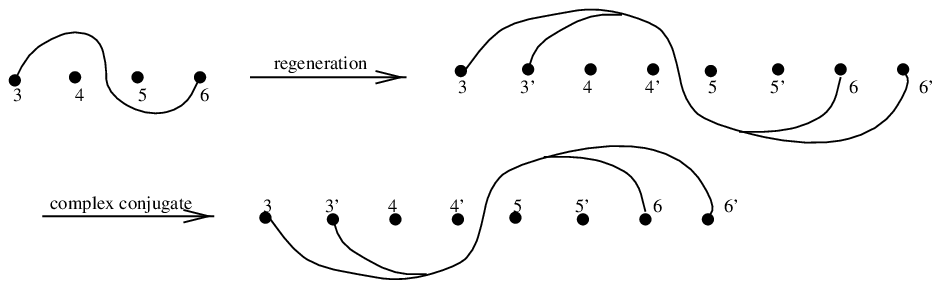}}}}$

$[\G_{4'}\G_{4}\G_{\check{3}3'} \G_4 \G_{4'}, \G_{66'}] = 1$ and
$[\G_{33'},\G_5 \G_{5'} \G_{6 \hat{6}} \G_{5'} \G_5]=1$.
\vspace{-1cm}

\item $\tilpat{Z}_{44',66'}^2$.
$\vcenter{\vspace{1cm}\hbox{\epsfxsize=9.75cm\epsfbox{\figs{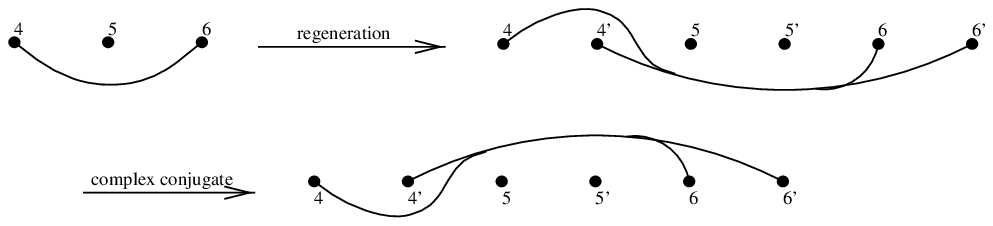}}}}$

$[\G_{\check{4}4'}, \G_{66'}] = 1$ and $[\G_{44'},\G_5 \G_{5'}
\G_{6 \hat{6}} \G_{5'} \G_5]=1$.
\end{itemize}
\end{tab}

\subsection{Checking Degrees}\label{degree}

Having computed the $C_i'$ (Subsection \ref{ssec:32}) and the
$H_{V_i}$ (Figure \ref{Vtab}), we obtain  a regenerated braid
monodromy factorization $\Dl_S^2=\prod\limiits_{i=1}^6 C'_i H_{V_i}$. To
verify that no factors are missing we compare degrees.  First,
since $S$ is a curve of degree $12$ (double the $6$ lines in
$S_0$), the braid $\Dl_S^2$ has degree $12\cdot 11=132$.  The six
monodromies $H_{V_i}$ each consist of three cusps and one branch point
for a combined degree of $6\cdot(3\cdot 3+1)=60$.  The $C'_i$
consist of four nodes for each parasitic intersection.  The nine
parasitic intersections (Table \ref{Ctab}) give a combined degree
of $9\cdot(4\cdot 2) = 72$. So together $\prod C'_iH_{V_i}$ has
also a degree of $60+72=132$, which proves that no factor was left
out.

\section{Invariance Theorems and $\Pitil_1$}\label{sec:4}

\subsection{The Invariance Theorem}\label{ssec:41}

Invariance properties are results concerning the behavior of a braid monodromy
factorization under conjugation by certain  elements of the braid group.
A factorization $g = g_1 \cdots g_k$ is said to be invariant under $h$ if
$g = g_1 \cdots g_k$ is Hurwitz equivalent to $(g_1)_h \cdots (g_k)_h$.
Geometrically this means that if a braid monodromy factorization of
$\Delta^2_{C}$ coming from a curve $S$ is invariant under $h$, then the
conjugate factorization is also a valid braid monodromy factorization for
$S$.

 The following rules \cite[Section 3]{MoTe4} give invariance
properties of commonly occurring subsets of braid monodromy
factorizations.  Factors of the third type do not appear in our
factorization.

(a) $\pat{Z}^2_{ii',jj'}$ is invariant under $\pat{Z}^p_{ii'}$ and $\pat{Z}^p_{jj'}$, $\forall p \in \Z$.

(b) $\pat{Z}^3_{i,jj'}$ is invariant under $\pat{Z}^p_{jj'}$, $\forall p \in \Z$.

(c) $\pat{Z}^1_{ij}$ is invariant under $\pat{Z}^p_{ii'} \pat{Z}^p_{jj'}$, $\forall p \in \Z$.

\begin{remark}
The elements $\pat{Z}_{ii'}$ and $\pat{Z}_{jj'}$ commute for all
$1\le i,j\le 6$ since the path from $i$ to $i'$ does not intersect
the path from $j$ to $j'$.
\end{remark}

\begin{theorem}[Invariance Theorem]\label{th2}
The braid monodromy factorization\break $\Delta^2_{12}= \prod\limiits^6_{i=1} C'_i H_{V_i}$
is invariant under  $\prod\limiits^6_{j=1} \pat{Z}^{m_j}_{jj'}$, for all $m_j \in \Z$.
\end{theorem}

\begin{proof}
It is sufficient to show that the $C'_i$ and the $H_{V_i}$ are
invariant individually.  Corollary 14 of \cite{MoTe4} proves that
each  $H_{V_i}$ is invariant under $\pat{Z}_{jj'}$, $1\le j\le 6$. Since
the $\pat{Z}_{jj'}$ all commute the invariance extends to arbitrary
products $\prod\limiits^6_{j=1} \pat{Z}^{m_j}_{jj'}$. The $C'_i$ are composed
of quadruples of factors $\tilpat{Z}_{kk',\ell\ell'}$, one from each
parasitic intersection. Lemma 16 of that paper shows that each
$\tilpat{Z}_{kk',\ell\ell'}$ is invariant under $\pat{Z}_{jj'}$, $1\le
j\le 6$. As before the invariance extends to products
$\prod\limiits^6_{j=1} \pat{Z}^{m_j}_{jj'}$. So the factorization
$\Delta^2_{12} = \prod\limiits^6_{i=1} C'_i H_{V_i}$ is invariant under
conjugation by these elements.
\end{proof}

We use $\G_{(j)}$ to denote any element of the set
$\{(\G_j)_{\pat{Z}_{jj'}^m}\}_{m \in \Z}$.  These elements are odd
length alternating products of $\G_j$ and $\G_{j'}$. Thus
$\G_{(j)}$ represents any element of
$\set{(\G_j\G_{j'})^p\G_j}_{p\in\Z}$. The original generators
$\G_j$ and $\G_{j'}$ are easily seen to be members of this set for
$p=0,-1$.

As an immediate consequence of the Invariance Theorem, any relation
satisfied by $\G_j$ is satisfied by any element of $\G_{(j)}$. 
This infinitely expands our collection of known relations in
$\Pitil_1$, however all of the new relations are consequences of
our original finite set of relations.
    $\Delta^2_{C}$ is also invariant under complex conjugation
\cite[Lemma 19]{MoTe4}, so we can use the complex conjugates
$\bar{H}_{V_i}$ and $\bar{C}'_i$ to derive additional relations.
Once again these relations are already implied by the existing
relations. On the other hand, many of the complex conjugate braids
in Table \ref{Ctab} have simpler paths than their counterparts so
they are a useful tool.

The paths corresponding
to the $H_{V_i}$ are already quite simple (see Figure \ref{Vtab})
so nothing is gained there by using
complex conjugates.

\subsection{A presentation for $\tilde\Pi_1$}\label{ssec:43}

Let $S$ be the regenerated branch curve and let $\pi_1 (\C^2 -
S,u_0)$ be the fundamental group of its complement in $\C^2$. We
know that this group is generated by the elements $\set{\G_j,
\G_{j'}}_{j=1}^6$.  Recall (Equation \eq{Pitildef}) that
$\Pitil_1 = {\pi_1(\C^2 - S,u_0) / \sg{\G^2_j, \G^2_{j'}}}$.

We have listed the braids $C_i$ (Table \ref{Ctab}) and $H_{V_i}$
(Figure \ref{Vtab}). These are the only braids in the factoization
of $\Delta_{C}^2$, as explained in Subsection \ref{degree}.

To the path of each braid
there correspond two elements of $\pi_1(\C^2-S,u_0)$,
as explained in Subsection \ref{ssec:32}.
{}From these, the van Kampen Theorem \cite{vK} produce the defining relations
of $\pi_1(\C^2-S,u_0)$.

\begin{theorem}\label{pres}
The group $\Pitil_1$ is generated by $\set{\G_j, \G_{j'}}_{j=1}^6$ with
the following relations:
\begin{eqnarray}
\Gamma_j^2 & = & 1 \qquad j=1,\dots,6 \label{sq}\\
\Gamma_{j'}^2 & = & 1 \qquad j=1,\dots,6 \label{sq2}\\
\G_{3'} & = & \G_1\G_{1'}\G_3\G_{1'}\G_1 \label{conjeq3t}\\
\G_{5'} & = & \G_4\G_{4'}\G_5\G_{4'}\G_4 \label{conjeq5t}\\
\G_{2} & = & \G_{6'}\G_{6}\G_{2'}\G_{6}\G_{6'} \label{conjeq2}\\
\G_{2'} & = & \G_1\G_{1'}\G_2\G_{1'}\G_1 \label{conjeq2t}\\
\G_{3} & = & \G_{4'}\G_{4}\G_{3'}\G_{4}\G_{4'} \label{conjeq3}\\
\G_{5} & = & \G_{6'}\G_{6}\G_{5'}\G_{6}\G_{6'} \label{conjeq5}\\
{}[\G_{ii'},\G_{jj'}] & = & 1 \qquad \mbox{if the lines $i,j$ are
disjoint in $X_0$} \label{commuteeq}\\
{}\G_{(i)}\G_{(j)}\G_{(i)}  & = & \G_{(j)}\G_{(i)}\G_{(j)} \qquad
\mbox{if the lines $i,j$ intersect}\label{tripeq}.
\end{eqnarray}
The enumeration of the lines is given in Figure \ref{sixe}.
$\Gamma_{ii'}$ represents either $\Gamma_i$ or $\Gamma_{i'}$, and
$\Gamma_{(i)}$ stands for any odd length word in the
infinite dihedral group $\sg{\Gamma_{i},\Gamma_{i'}}$.
\end{theorem}
\begin{proof}
Relations \eqs{sq}{sq2} hold in $\Pitil_1$ by assumption. The
other relations hold in $\pi_1(\C^2-S,u_0)$. To see this, we now
list the relations induced by the $H_{V_i}$ (Table \ref{HV}). Recall that each
$H_{V_i}$ is a product of regenerated braids induced from one
branch point (the second path in each part of Figure \ref{Vtab}),
and three cusps (condensed in the first path of each part). 
Applying the van Kampen Theorem \cite{vK}, we have two types of
relations. The relations \eqs{conjeq3t}{conjeq5} are derived from
the branch point braids, and the triple relations $\G_i\G_j\G_i =
\G_j\G_i\G_j$ from the three other braids. Using the Invariance
Theorem \ref{th2} to expand the pathes, we get \eq{tripeq} in its
full generality. It remains to prove Equation~\eq{commuteeq}.
Note that the relations 
$[\G_{ii'},\G_{jj'}] = 1$ and $[\G_{(i)},\G_{(j)}] = 1$ imply each other.

We now consider the complex conjugates in Table \ref{Ctab}. The
relations $[\G_{(2)},\G_{(3)}]$, $[\G_{(2)},\G_{(4)}]$,
$[\G_{(1)},\G_{(5)}]$, $[\G_{(2)},\G_{(5)}]$,
$[\G_{(3)},\G_{(5)}]$,  and $[\G_{(4)},\G_{(6)}]$
 appear rather
directly. We must derive the other commutators, namely
$[\G_{(1)},\G_{(4)}]$, $[\G_{(1)},\G_{(6)}]$, and
$[\G_{(3)},\G_{(6)}]$.

By the second part of Table \ref{Ctab}, 
$$[\G_{2'}\G_2\G_{(1)}\G_2\G_{2'},\G_{(4)}]=
\G_{2'}\G_2\G_{(1)}\G_2\G_{2'}\G_{(4)}\G_{2'}\G_2\G_{(1)}\G_2\G_{2'}\G_{(4)},$$
but since $[\G_{(2)},\G_{(4)}] = 1$ we get
$$\G_{2'}\G_2\G_{(1)}\G_{(4)}\G_2\G_{2'}\G_{2'}\G_2\G_{(1)}\G_{(4)}\G_2\G_{2'}=
\G_{2'}\G_2\G_{(1)}\G_{(4)}\G_{(1)}\G_{(4)}\G_2\G_{2'}$$ from
which $[\G_{(1)}, \G_{(4)}] = 1$ follows.

Now, by the seventh part of the table, 
$$[\G_{(1)},\G_5\G_{5'}\G_{(6)}\G_{5'}\G_5]=
\G_{(1)}\G_5\G_{5'}\G_{(6)}\G_{5'}
\G_5\G_{(1)}\G_5\G_{5'}\G_{(6)}\G_{5'}\G_5$$ but since
$[\G_{(1)},\G_{(5)}] = 1$ we get
$$\G_5\G_{5'}\G_{(1)}\G_{(6)}\G_{5'}\G_5\G_5\G_{5'}\G_{(1)}\G_{(6)}\G_{5'}\G_5=
\G_5\G_{5'}\G_{(1)}\G_{(6)}\G_{(1)}\G_{(6)}\G_5\G_{5'}$$ from
which we get $[\G_{(1)}, \G_{(6)}] = 1$.

In the same way, since $\G_{(3)}$ and $\G_{(5)}$ commute we can
get $[\G_{(3)}, \G_{(6)}] = 1$ from the relation
$[\G_{(3)},\G_5\G_{5'}\G_{(6)}\G_{5'}\G_5]$ (the eighth part of the table).
This finishes the proof of \eq{commuteeq}.
\end{proof}

\section{The homomorphism $\psi$}\label{sec:5}

$X^{\Aff}-f^{-1}(S)$ is a degree $6$ covering space of $\C^2-S$.
Let
$$\psi\colon\pi_1(\C^2-S,u_0)\to S_6$$ be the permutation monodromy of this
cover. As before let $\pi\colon\C^2\to\C$ be a generic projection
and choose $u\in\C$ such that $S$ is unramified over $u$. For
surfaces $X$ close to the degenerated $X_0$ the points $i$ and
$i'$ will be close to each other in $\C_u$. Finally choose a
point $u_0 \in \C_u$.

We wish to determine what happens to the six preimages $f^{-1}(u_0)$
in $X$ as they follow the lifts of $\G_i$ and $\G_{i'}$.  Again, for
surfaces $X$ close to the degenerate $X_0$ these six points inherit a
unique
enumeration based on which numbered plane of $X_0$ they are nearest.
This enumeration remains valid all along $\alpha_i$ so we need only
consider the monodromy around $\eta_i$ and $\eta'_i$.  Take a small
neighborhood $U_i\sbt\C_u$ of $i$ and $i'$.
We can reduce the dimension of the question by restricting to
$f^{-1}(U_i)$ which is a branched cover of $U_i$, branched over
$i$ and $i'$.  It is clear that $f_0^{-1}(U_i)$ has a
simple node over $i$ where two planes containing $\hat L_i$
meet.  When $i$ divides, the node will factor into two simple branch
points over $i$ and
$i'$ involving the same sheets which met at the node.  Thus we see that
if $P_k$ and $P_\ell$ intersect in $\hat L_i$ then
$\psi(\G_i)=\psi(\G_{i'})=(k\ \ell)$. Specifically, $\psi$ is defined by
\begin{defn}\label{defpsi}
The map $\psi\colon\pi_1(\C^2-S,u_0)\to S_6$ is given by
\begin{eqnarray*}
\psi(\G_1)=\psi(\G_{1'})& = & (13),\\
\psi(\G_2)=\psi(\G_{2'})& = &(15),\\
\psi(\G_3)=\psi(\G_{3'})&=&(23),\\
\psi(\G_4)=\psi(\G_{4'})&=&(26),\\
\psi(\G_5)=\psi(\G_{5'})&=&(46),\\
\psi(\G_6)=\psi(\G_{6'})&=&(45).
\end{eqnarray*}
\end{defn}

Figure \ref{sixonceagain} depicts
the simplicial complex of $X_0$ with the planes
and intersection lines numbered. {}From this we can determine the
values of $\psi$ on the generators $\G_i$ and $\G_{i'}$. 
Figure~\ref{hex} below
gives another graphical representation for
$\psi$, in which $\Gamma_i$ connects the two vertices $\alpha,\beta$ 
defined by $\psi(\Gamma_i) = (\alpha\beta)$.

{\nocolon\FIGUREx{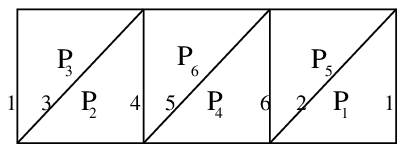}{sixonceagain}{5cm}}

The reader may wish to check that $\psi$ is well defined (testing
the relations given in Theorem \ref{pres}), but this is of course
guaranteed by the theory. {}From the definition $\psi$ is clearly
surjective.

Since $\psi(\G^2_j) = 1$ and $\psi(\G^2_{j'}) = 1$, $\psi$ also defines
a map $\Pitil_1\to S_6$.  We will use $\psi$ to denote this map as
well. Let  $\mathcal{A}$ be the kernel of $\psi\colon\Pitil_1
\rightarrow S_6$.  We have a short exact sequence sequence
\begin{equation}\label{short}
1 \longrightarrow \mathcal{A}
 \longrightarrow \Pitil_1 \stackrel{\psi}{\longrightarrow}  S_6 \rightarrow 1.
 \end{equation}
\begin{theorem}{\rm\cite{MoTe1}}\qua
The fundamental group $\pi_1(X^{\Aff}_{\Gal})$ is isomorphic to
$\mathcal{A}$,\break where $X^{\Aff}_{\Gal}$ is the affine part of 
$X_{\Gal}$ the
Galois cover of $X$ with respect to a generic projection onto
$\C\P^2$.
\end{theorem}

\section{A Coxeter subgroup of $\GG$}\label{sec:6} 

Our next step in identifying the group $\GG$ is to study a natural
subgroup, which happens to be a Coxeter group.

Let $\theta \co \GG\rightarrow C$ be the map defined by $\theta(\G_j) = \theta(\G_{j'}) = \gamma_j$.
The resulting group $C =  \Im(\theta)$ is formally defined by the generators $\gamma_1,\dots,\gamma_6$ and
the relations obtained by applying $\theta$ to the relations of $\GG$.

Since we have $\psi(\G_j) = \psi(\G_{j'})$, $\psi$ splits through
$\theta$:
defining $\rho \co C \rightarrow S_6$ by 
$\rho(\gamma_j) = \psi(\G_j)$, we have that $\psi = \rho \circ \theta$.

\begin{lemma}\label{lm4}
In terms of the generators $\g_j$, $C$ has the following presentation
\begin{eqnarray*}
C & = & \la \g_1, \cdots ,\g_6\ |\ \g^2_i = 1,\\
&& \Trip{\g_1}{\g_3},\ \Trip{\g_3}{\g_4},\ \Trip{\g_4}{\g_5},
\nonumber \\
&& \Trip{\g_5}{\g_6},\ \Trip{\g_6}{\g_2},\ \Trip{\g_2}{\g_1},
\\ \nonumber && [\g_1, \g_4] , [\g_1, \g_5], [\g_1, \g_6], [\g_3,
\g_5], [\g_3, \g_6], [\g_3, \g_2], [\g_4, \g_6],[\g_4,
\g_2],[\g_5, \g_2] \ra    .
\end{eqnarray*}
\end{lemma}
\begin{proof}
We only need to apply $\theta$ on the relations of $\Pitil_1$
given in Theorem \ref{pres}. Each of the relations in $\Pitil_1$
descends to a relation in $C$. The branch points all give
relations of the form $\G_{3'} =
\G_{1}\G_{1'}\G_{3}\G_{1'}\G_{1}$. When we equate $\G_j = \G_{j'}$
these relations vanish. If $L_i,L_j$ intersect, then the relations 
$\Trip{\G_{(i)}}{\G_{(j)}}$ descend to 
$\Trip{\g_i}{\g_j}$. Similarly if the lines $\hat L_i,\hat L_j$ 
are disjoint, we get $[\g_i, \g_j]$.
\end{proof}

This presentation is easier to remember using Figure~\ref{hex}: $\g_i,\g_j$ 
satisfy the triple relation if the corresponding lines intersect in a common
vertex, and commute otherwise.

\begin{figure}[ht!]\small
\begin{equation*}
\xymatrix{ 
	{}
	& 3 \ar@{-}[r]^{\gamma_3}
	& 2 \ar@{-}[rd]^{\gamma_4}
	& {}
\\
	1 \ar@{-}[ru]^{\gamma_1}
	& {}
	& {}
	& 6 \ar@{-}[ld]^{\gamma_5}
\\
	{}
	& 5 \ar@{-}[lu]^{\gamma_2}
	& 4 \ar@{-}[l]^{\gamma_6}
	& {}
}
\end{equation*}\vspace{-6mm}
\nocolon\caption{}\label{hex}
\end{figure}

We will use Reidemeister-Schreier method to compute $\mathcal
A$. For this we need  to split $\theta$.

\begin{lemma}\label{lm5}
$\theta$ is split by the map $s\colon C\to\Pitil_1$ defined by
$\g_j \longmapsto \G_j$.
\end{lemma}
\begin{proof}
{}From the definition of $\theta$ it is clear that $\theta \circ
s$ is the identity on $C$, so it remains to check that $s$ is well
defined. In Lemma \ref{lm4} we gave a presentation of $C$.  To
prove that $\g_j \longmapsto \G_j$ is a splitting we must show
that $\{\G_j\}$ also satisfies the relations in $\Pitil_1$. $s$
respects $\g^2_j=1$, since $\G^2_j = 1$ in $\Pitil_1$. Let $i,j$ be indices
\st\ $L_i,L_j$ intersect, then $\Trip{\g_i}{\g_j}$ is respected
since  in $\Pitil_1$ we have $\Trip{\G_{(i)}}{\G_{(j)}}$ so
specifically $\Trip{\G_{i}}{\G_{j}}$. Finally if $i,j$ are  \st\
$\hat L_i \cap \hat L_j = \phi$,  $[\g_i, \g_j]$ is respected since
$[\G_{ii'},\G_{jj'}] = 1$ for disjoint $i, j$.
\end{proof}

Observe that Lemma \ref{lm4} presents $C$ as a Coxeter group on the generators
$\gamma_1,\gamma_3,\gamma_4,\gamma_5,\gamma_6,\gamma_2$, with a hexagon 
(the dual of that shown in Figure \ref{hex}) as the Coxeter-Dynkin
diagram of the group. The fundamental group of this defining graph
is of course $\Z$.
In previous works on the fundamental groups of Galois covers, 
the group $C$ defined in a similar
manner to what we define here, always happen to be equal to the symmetric group
$S_n$ (where $n$ is the number of planes in the degeneration). Here,
the map from $C$ to $S_6$ is certainly not injective 
($C$ is known \cite{Bourbaki} to be the group $S_6 \semidirect \Z^5$,
with an action of $S_6$ on $\Z^5$ by the nontrivial component of the standard
representation). The connection of this fact to the fundamental group of $S_0$
is explained in more details in \cite{TxT}.

It will be useful for us to have a concrete isomorphism of $C$ and $S_6 \semidirect \Z^5$.

\begin{lemma}\label{Cstructure}
$C \cong S_6 \semidirect \Z^5$ where $\Z^5$ is the nontrivial
component of the standard representation.
\end{lemma}
\begin{proof}
First note that $\sg{\gamma_2,\dots,\gamma_6}$ is the parabolic
subgroup of $C$ corresponding to the Dynkin diagram of type $A_5$,
so that $C_0 = \sg{\gamma_2,\dots,\gamma_6} = S_6$. We will
therefore identify the subgroup $C_0$ with the symmetric group
$S_6$ (using $\psi$ as the identifying map). Next, note that
$\psi(\g_1) = (13)$, so we set $x = (13)\g_1$, and consider the
presentation of $C$ on the new set of generators, namely
$x,\g_2,\dots,\g_6$. Substituting $\g_1 = (13)x$ in the
presentation of Lemma \ref{lm4}, we obtain $C = \sg{x,S_6}$,
with the relations%
\begin{eqnarray*}
(13)x(23)(13)x & = & (23)(13)x(23), \\
(13)x(15)(13)x & = & (15)(13)x(15), \\
(26)x(26) & = & x, \\
(46)x(46) & = & x, \\
(45)x(45) & = & x.
\end{eqnarray*}
Define $x_\s = \s x \s^{-1}$, then the fact that $x$ commute with
$\sg{(26),(46),(45)} = S_\set{2,4,5,6}$ shows that $x_\s$ actually
depends only on $\s^{-1}(1),\s^{-1}(3)$. We can thus define
$x_{k\ell} = x_{\s}$ for some $\s \in \sg{\gamma_2,\dots,\gamma_6}$
such that $\s(k)=1$ and
$\s(\ell) = 3$ (so in particular $x_{13} = x$). With this
definition one checks that $\s^{-1}x_{k\ell}\s =
x_{\s(k),\s(\ell)}$. Adding this last relation as a definition of
the $x_{k\ell}$, we obtain the following presentation: $C =
\sg{x_{k\ell}, S_6}$, with the
relations%
\begin{eqnarray*}
(13)x_{13}(23)(13)x_{13} & = & (23)(13)x_{13}(23), \\
(13)x_{13}(15)(13)x_{13} & = & (15)(13)x_{13}(15), \\
\s^{-1}x_{k\ell}\s & = & x_{\s(k),\s(\ell)}.
\end{eqnarray*}
Now, the first two relations translate to%
\begin{eqnarray*}
x_{32}x_{13} & = & x_{12}, \\
x_{51}x_{13} & = & x_{53},
\end{eqnarray*}
which after conjugating by an arbitrary $\s$ give%
\begin{eqnarray*}
x_{ij}x_{ki} & = & x_{kj}, \\
x_{ki}x_{ij} & = & x_{kj},
\end{eqnarray*}
which shows that $\sg{x_{k\ell}}$ is generated by
$x_{12},x_{13},\dots,x_{16}$ and is commutative (using the fact that the 
$x_{1i}$ commute).
Thus $\sg{x_{k\ell}} = \Z^5$, and $C
= \sg{\Z^5,S_6}$ is the asserted group.
\end{proof}

 The inclusion $S_6 \hookrightarrow C$ defined by sending
the transpositions $(15)$, $(23)$, $(26)$, $(46)$ and $(45)$ to $\g_2,\dots,
\g_6$ respectively, splits the projection $\rho \co C \rightarrow
S_6$. {}From now on we identify $S_6$ with the subgroup
$\sg{\g_2,\dots,\g_6}$ of $C$, as well as the subgroup
$\sg{\G_2,\dots,\G_6}$ of $\GG$.

\begin{corollary}\label{sesissplit}
The sequence \eq{short} is split (by the composition of the maps
$S_6 \hookrightarrow C$ and $s \co C \hookrightarrow \GG$). We denote
the splitting map by $\varphi$.
\end{corollary}

\section{The kernel of $\psi$}\label{sec:7}

We use the Reidemeister-Schreier method to find a presentation for
the kernel $\mathcal{A}$ of the map $\psi\colon\Pitil_1
\rightarrow S_6$. 
Let $L = \Ker(\theta \co \GG \rightarrow C)$ and $K = \Ker(\rho co C \rightarrow S_6)$, and consider the diagram of Figure \ref{ninediag}, in which the rows 
are exact by definition of $\mathcal{A}$ and $K$, and the middle colomn by definition of
$L$. The equality $L = L$ in the diagram follows from the nine lemma.
Then, since $1 \rightarrow L \rightarrow \GG \rightarrow C
\rightarrow 1$ splits (Lemma \ref{lm5}), we have that
$\mathcal{A}$ is a semidirect product of $L = \Ker(\theta \co \GG
\rightarrow C)$ and $K = \Ker(\rho \co C \rightarrow S_6)$,
which is isomorphic to $\Z^5$ by Lemma \ref{Cstructure}.

\begin{figure}[ht!]\small
\begin{equation*}
\xymatrix{ 
	{}
	& 1 \ar@{->}[d]
	& 1 \ar@{->}[d]
	& {}
	& {}
\\
	{}
	& L \ar@{=}[r] \ar@{->}[d]
	& L \ar@{->}[d]
	& {}
	& {}
\\
	{} 1 \ar@{->}[r]
	& \mathcal{A} \ar@{->}[r] \ar@<.5ex>@{->}[d]
	& \Pitil_1 \ar@<.5ex>@{->}[r]^{\psi} \ar@<.5ex>@{->}[d]^{\theta}
	& S_6 \ar@{->}[r] \ar@<.5ex>@{->}[l]^{\varphi}
	& 1
\\
	{} 1 \ar@{->}[r]
	& K \ar@{->}[d] \ar@<.5ex>@{->}[u] \ar@{->}[r]
	& C \ar@{->}[d] \ar@<.5ex>@{->}[r]^{\rho} \ar@<.5ex>@{->}[u]
	& S_6 \ar@{->}[r] \ar@<.5ex>@{->}[l]
	& 1
\\
	{}
	& 1
	& 1
	& {}
	& {}
}
\end{equation*}\vspace{-5mm}
\nocolon\caption{}\label{ninediag}
\end{figure}

\subsection{The Reidemeister-Schreier method}\label{ss:RMS}

Let
$$ 1 \rightarrow K \rightarrow G \stackrel{\phi}{\rightarrow} H \rightarrow 1$$
be a short exact sequence, split by $\rho \co H \rightarrow G$.
Assume that $G$ is finitely generated, with generators
$a_1, \cdots , a_n$.
Then $\rho \phi (g)$ is a representative for $g \in G$ in its class modulo $H$.
It is easy to see that $K$ is generated by the elements
$g a_i{(\rho \phi(g a_i))}^{-1}$, $1 \leq i \leq n$, $g \in \rho(H)$.
Now,
$ ga_i{(\rho \phi (ga_i))}^{-1}= ga_i(\rho \phi ({a_i}))^{-1}
{(\rho \phi (g))}^{-1}=
 ga_i (\rho \phi ({a_i}))^{-1} g^{-1}$, because
$\phi \rho$ is the identity on $H$.
For $g \in \rho(H)$, denote the generators above by
$$\g(g, a_i) = g a_i (\rho \phi ({a_i}))^{-1} g^{-1}.$$

The relations of $G$ can be translated into expressions in these 
generators by the following process.
If the word
$\omega = a_{i_1} \cdots a_{i_t}$ represents an element of $K$
then $\omega$ can be rewritten as the product
$$\tau (\omega) = \g(1,a_{i_1}) \g(\rho\phi(a_{i_1}),a_2)
\cdots
\g(\rho \phi(a_{i_1} \cdots a_{i_{t-1}}),  a_{i_t}).$$

\begin{theorem}[Reidmeister-Schreier]\label{RMS}
Let $\set{R}$ be a complete set of relations for $G$.
Then
$K = \Ker(\phi)$
is generated by the $\g(g,a_i)$ ($1 \leq i \leq n$, $g \in H$), with
the relations
$\{\tau(trt^{-1})\}_{r \in R, t \in \rho(H)}$.
\end{theorem}

We will use this method to investigate $L = \ker \theta$ and
$\mathcal{A} = \ker \psi$.

\subsection{Generators for $L = \Ker \theta$}\label{ss:presL}
\def\bG{\bar{\Gamma}}

For $c \in \Im(s) = \sg{\G_1,\dots,\G_6}$, we
let
\begin{equation}\label{defAcj}
A_{c,j} = c \G_{j} \G_{j'} c^{-1}.
\end{equation}

We start with the following:
\begin{corollary}\label{genL}
The group $L = \Ker(\theta)$ is generated by
$\{A_{c,j}\}_{1\le j\le6,\, c\in C}$.
\end{corollary}
\begin{proof}
By Theorem \ref{RMS}, $L$ is generated by elements of the form\nl
$c\G_j {(s \theta({\G}_j))}^{-1} c^{-1}$ and $c \G_{j'} (s \theta
({\G}^{-1}_{j'})) c^{-1}$ for all $1 \leq  j \leq 6$ and $c \in
C$. We compute $s \theta(\G_{j}) = s(\g_j) = \G_j$ and $ s
\theta(\G_{j'}) = s(\g_j) = \G_j$, so the generators are $c \G_j
\G^{-1}_j c^{-1} = I$ and $c \G_{j'} \G^{-1}_j c^{-1} = c \G_{j'}
\G_j c^{-1}$ using $\G^2_j = I$.
\end{proof}

This set of generators is
highly redundant as we shall later see, but for now we turn our
attention to $\mathcal{A} = \Ker \psi$.

\subsection{Generators for $\mathcal{A} = \Ker \psi$}\label{sec:10}

By Theorem \ref{RMS}, $\mathcal{A}$ is generated
by the elements $\s \G_j (\varphi\psi(\G_j))^{-1}
 \s^{-1} $ and\nl $\s \G_{j'}
 (\varphi\psi(\G_{j'}))^{-1} \s^{-1}$, $1 \leq j \leq 6$, $\s \in
S_6$.  Again we compute $\varphi\psi(\G_j)$ and $\varphi
\psi(\G_{j'})$. Recall that $\sg{\G_2,\dots,\G_6} = S_6$ is the
image of $\varphi$ (Corollary \ref{sesissplit}), so for $j \neq 1$
we get $\varphi\psi(\G_j) = \varphi \psi(\G_{j'}) = \G_j$, and the
generators are
\begin{equation}\label{defAsj}
A_{\s,j}=\s \G_{j}\G_{j'}\s^{-1}.
\end{equation}
This agrees with our previous definition of $A_{c,j}$ for $c \in S_6$, see
Equation \eq{defAcj}.

The permutation $(1\,3)$ can be expressed in terms of the generators
of $S_6$ corresponding to $\Gamma_2,\dots,\Gamma_6$ as follows:
$$(1\,3) = (1\,5)(5\,4)(4\,6)(2\,6)(2\,3)(2\,6)(4\,6)(5\,4)(1\,5),$$
so for $j = 1$ we have that
$$\varphi\psi(\G_1) = 
\varphi((15)(54)(46)(62)(23)\cdots(15)) = \G_2\G_6\G_5\G_4
\G_3\G_4\G_5\G_6\G_2.$$ Likewise,
$\varphi\psi(\G_1)=\varphi\psi(\G_{1'})$ since $\psi(\G_j) =
\psi(\G_{j'})$. So we get generators
\begin{eqnarray}
X_{\s} & = & \s\G_2\G_6\G_5\G_4\G_3\G_4\G_5\G_6\G_2\G_1\s^{-1}\label{Xdef}\\
B_{\s} & = & \s\G_2\G_6\G_5\G_4\G_3\G_4\G_5\G_6\G_2\G_{1'}\s^{-1}.\label{Bdef}
\end{eqnarray}

Since $X_\s^{-1}B_\s =
\s \G_1 \G_{1'} \s^{-1}  = A_{\s,1}$,
we have the following result:
\begin{corollary}
The group $\mathcal{A} = \Ker(\psi)$ is generated by
$A_{\s,j}$, $X_\s$, for $\s\in S_6$ and $j = 1,\dots,6$.
\end{corollary}

Notice that we are now conjugating only by permutations $\s\in S_6$
instead of all elements $c\in C$ (as in Corollary \ref{genL}) so this is a finite set of generators.

\section{A better set of generators for $\mathcal{A}$}\label{sec:8}

We first show that $A_{\s,j}$ are not needed for $j = 2,\dots,6$

\begin{theorem}\label{th9}
$\mathcal{A}$ is generated by $\{A_{\s,1}, X_\s\}$.
\end{theorem}
\begin{proof}
This follows immediately from the relations proven below.
\end{proof}

\begin{tab}\label{Vtab-rel} 
We have the following relations:
\begin{eqnarray}
A_{\s,3} & = & A_{\s(23),1}A_{\s,1}^{-1} \label{AA1}\\
A_{\s,3} & = & A_{\s(26)(23),4} \\
A_{\s,5} & = & A_{\s(26)(46),4} \\
A_{\s,5} & = & A_{\s(45)(46),6} \\
A_{\s,2} & = & A_{\s(45)(15),6} \\
A_{\s,2} & = & A_{\s(15),1}A_{\s,1}^{-1}.\label{AA6}
\end{eqnarray}
\end{tab}
\begin{proof}
We use the relations of Theorem \ref{pres}. Let $I$ denote the
identity element of $S_6$, so that by definition $A_{I,j} =
\Gamma_j\Gamma_{j'}$. %
{}From \eq{conjeq3t} 
we have $1 = \G_{3'}\G_1\G_{1'}\G_3\G_{1'}\G_1=
(\G_{3'}\G_3)(\G_3\G_1\G_{1'}\G_3)(\G_{1'}\G_1)=
A_{I,3}^{-1}A_{(23),1}A_{I,1}^{-1}$ so we get
$A_{I,3}$\break$=A_{(23),1}A_{I,1}^{-1}$.

{}From \eq{conjeq5t} we have $1
 = \G_{5'}\G_4\G_{4'}\G_5\G_{4'}\G_4 =
(\G_{5'}\G_5) \G_5\G_4\G_{4'}\G_5\G_{4'}\G_4$\nl$=
(\G_{5'}\G_5) \G_5\G_4\G_{5}\G_{4'}\G_{5}\G_4 =
(\G_{5'}\G_5) \G_4\G_5\G_{4}\G_{4'}\G_{5}\G_4 =
A_{I,5}^{-1} A_{(26)(46),4}$.

{}From \eq{conjeq2} 
we have $1=\G_2\G_{6'}\G_6\G_{2'}\G_6\G_{6'}=
(\G_2\G_{6'}\G_6\G_2)(\G_2\G_{2'})(\G_6\G_{6'})$
so we get $A_{I,2} = \G_2\G_{2'} = \G_2\G_6\G_{6'}\G_2\G_{6'}\G_6 = 
\G_2\G_6\G_{2}\G_{6'}\G_{2}\G_6 = 
\G_6\G_2\G_{6}\G_{6'}\G_{2}\G_6 = A_{(45)(15),6}$.

{}From \eq{conjeq2t} 
we have $1=\G_{2'}\G_1\G_{1'}\G_2\G_{1'}\G_1=
(\G_{2'}\G_2)(\G_2\G_1\G_{1'}\G_2)(\G_{1'}\G_1)$\nl$=
A_{I,2}^{-1}A_{(15),1}A_{I,1}^{-1}$ so we get
$A_{I,2}=A_{(15),1}A_{I,1}^{-1}$.

{}From \eq{conjeq3} 
we have $1=\G_3\G_{4'}\G_4\G_{3'}\G_4\G_{4'}=
(\G_3\G_{4'}\G_4\G_3)(\G_3\G_{3'})(\G_4\G_{4'})$
so we get $\G_3 \G_{3'} = \G_3\G_4\G_{4'}\G_3\G_{4'}\G_4 = 
\G_3\G_4\G_{3}\G_{4'}\G_{3}\G_4 = 
\G_4\G_3\G_{4}\G_{4'}\G_{3}\G_4 = A_{(26)(23),4}$.

Finally from \eq{conjeq5} 
we have $1=\G_5\G_{6'}\G_6\G_{5'}\G_6\G_{6'}=
(\G_5\G_{6'}\G_6\G_5)(\G_5\G_{5'})(\G_6\G_{6'})$
so we get 
$\G_5\G_{5'} = \G_5\G_{6}\G_{6'}\G_5\G_{6'}\G_{6}
= \G_5\G_{6}\G_{5}\G_{6'}\G_{5}\G_{6}
= \G_6\G_{5}\G_{6}\G_{6'}\G_{5}\G_{6} = A_{(45)(46),6}$.
\end{proof}

One may be tempted to use 
$(\G_3\G_1\G_{1'}\G_3)(\G_{1'}\G_1) = \G_1\G_3\G_{1}\G_{1'}\G_{3}\G_1$ 
in a similar manner to the other cases, to rewrite
$A_{\s (23),1}A_{\s,1}^{-1}$ of Equation \eq{AA1} as a single element of the
form $A_{\s,1}$; however note that in the definition \eq{defAsj} we require 
$\s \in S_6 = \sg{\G_2,\dots,\G_6}$, so we do not have
the equality $\G_1\G_3\G_{1}\G_{1'}\G_{3}\G_1 = A_{(23)(13),1}$. 
The same remark applies for  Equation \eq{AA6}.

Iterating the relations of Table \ref{Vtab-rel}, we obtain a new
relation for $\sg{A_{\s,1}}$:
\begin{eqnarray*}
A_{\s(23),1}A_{\s,1}^{-1}& = & A_{\s,3} \\
& = & A_{\s(26)(23),4} \\
& = & A_{\s(26)(23)(46)(26),5}
\\
& = & A_{\s(26)(23)(46)(26)(45)(46),6} \\
& = & A_{\s(26)(23)(46)(26)(45)(46)(15)(45),2} \\
& = & A_{\s(36)(24)(56)(14),2} \\
& = & A_{\s(36)(24)(56)(14)(15),1} A_{\s(36)(24)(56)(14),1}^{-1}
\end{eqnarray*}
which may be rewritten as
\begin{equation}\label{eqA1}
A_{\s(23),1}A_{\s,1}^{-1}=A_{\s(142563),1}A_{\s(142)(356),1}^{-1}
\end{equation}

\begin{lemma}\label{Xkllem}
For every $\s \in S_6$, $X_\s$ depends only on $\s^{-1}(1)$ and $\s^{-1}(3)$.
\end{lemma}
\begin{proof}
Viewing $S_6 \subset C$ as subgroups of $\Pitil_1$ (using the
embedding $s\co C \rightarrow \Pitil_1$), the elements $X_\s$ belong
to $C = \sg{\G_1,\G_2,\dots,\G_6}$ by their definition \eq{Xdef}.

Applying the isomorphism
$C \cong  S_6 \semidirect \Z^5$ of Lemma \ref{Cstructure}, we see that
\begin{eqnarray*}
X_{\s} & = & \s\G_2\G_6\G_5\G_4\G_3\G_4\G_5\G_6\G_2\G_1\s^{-1} \\
	& = & s(\s(15)(54)(46)(62)(23)(62)(46)(54)(15)(13)x_{13}\s^{-1}) \\
	& = & s(\s(13)(13)x_{13}\s^{-1}) \\
	& = & s(\s x_{13}\s^{-1}) 
	 =  s(x_{\s^{-1}(1)\s^{-1}(3)}).
\end{eqnarray*}
\end{proof}

A similar result holds for $\set{A_{\s,1}}$.

\begin{lemma}\label{Akllem}
For every $\s \in S_6$, $A_{\s,1}$
depends only on $\s^{-1}(1)$ and $\s^{-1}(3)$.
\end{lemma}
\begin{proof}
If $\s_1^{-1}(1) = \s_2^{-1}(1)$ and $\s_1^{-1}(3) = \s_2^{-1}(3)$,
then $\tau = \s_1^{-1} \s_2$ stabilizes $1,3$, so $\tau \in
S_\set{2,6,4,5} = \sg{\G_4,\G_5,\G_6}$ which commute with both
$\G_1$ and $\G_{1'}$. By definition \eq{defAsj}, $A_{\s_2,1} =
\s_2 \G_1 \G_{1'} \s_2^{-1} = \s_1 \tau \G_1 \G_{1'} \tau^{-1}
\s_1^{-1} = \s_1 \G_1 \G_{1'} \s_1^{-1} = A_{\s_1,1}$.
\end{proof}

We can thus define
\begin{defn}\label{AXkldef}
For $k,\ell = 1,\dots,6$, $A_{k\ell}$ and $X_{k\ell}$ are defined by
\begin{eqnarray}
X_{k\ell} & = & \s \G_2\G_6\G_5\G_4\G_3\G_4\G_5\G_6\G_2 \G_1 \s^{-1} \label{Xkl} \\
A_{k\ell} & = & \s \G_1\G_{1'} \s^{-1} \label{Akl}
\end{eqnarray}
where $\s \in S_6 = \sg{\G_2,\dots,\G_6}$ is any permutation such that
$\s(k)=1$ and $\s(\ell) = 3$.
\end{defn}

We need to know the action of $S_6$ on these generators:

\begin{proposition}
For every $\s \in S_6 = \sg{\G_2,\dots,\G_6}$ and $k,\ell = 1,\dots,6$,
we have that
\begin{eqnarray}
\s^{-1} A_{k\ell} \s & = & A_{\s(k),\s(\ell)}\label{As6act}\\
\s^{-1} X_{k\ell} \s & = & X_{\s(k),\s(\ell)}\label{Xs6act}
\end{eqnarray}
\end{proposition}
\begin{proof}
Let $\tau \in S_6$ be such that $\tau(k)=1$ and $\tau(\ell) =3$. Since
$A_{13} = \Gamma_1 \Gamma_{1'}$ by Definition \ref{AXkldef}, we have 
$A_{k\ell} = \tau A_{13} \tau^{-1}$
and $\s^{-1} A_{k \ell} \s = \s^{-1} \tau A_{13} \tau^{-1} \s = 
A_{\tau^{-1}\s(1),\tau^{-1}\s(3)} = A_{\s(k)\s(\ell)}$.
The same proof works for the $X_{k\ell}$.
\end{proof}

Note that $B_{k\ell} = X_{k \ell} A_{k \ell}$ can be defined in
the same manner, and have the same $S_6$-action.

{}From Theorem \ref{th9} (with a little help from Lemma \ref{Xkllem} and Lemma \ref{Akllem}), we obtain
\begin{corollary}\label{co12}
The group $\mathcal{A}$ is generated by
$\{A_{k\ell},X_{k\ell}\}_{1 \leq k,\ell \leq 6}$.
\end{corollary}

Since $\theta(X_{k\ell}) = x_{k\ell}$, we already
proved
\begin{corollary}\label{lm10}
The elements $\set{\theta(X_{k\ell})}$ generate $K =
\Ker(\rho\colon C\rightarrow S_6)$.
\end{corollary}

In the new language, Equation \eq{eqA1} (for $\s = I$) can be
written as $A_{12}A_{13}^{-1}=A_{36}A_{26}^{-1}$, so conjugating
we get
\begin{equation}\label{Akl:rel}
A_{ij}A_{ik}^{-1}=A_{k\ell}A_{j\ell}^{-1}
\end{equation}
for any four distinct indices $i,j,k,\ell$.  Using three consecutive
applications of the relation (\ref{Akl:rel}) we can also allow $i=\ell$,
and using just two applications we can get:
\begin{equation}\label{Akl2:rel}
A_{ij}A_{ik}^{-1}=A_{\ell j}A_{\ell k}^{-1}
\end{equation}
for any distinct indices $i,j,k$ and $\ell\ne j,k$.  In view of
(\ref{Akl:rel}) and (\ref{Akl2:rel}) and Table~\ref{Vtab-rel} we can
write a translation table for the remaining generators $A_{\s,j}$ for
$j\ne 1$.
\begin{tab}\label{Vtrans}
$A_{\s,j}$ in terms of $A_{k\ell}$
\begin{eqnarray}
A_{\s,1} & = & \s A_{13}\s^{-1} \\
A_{\s,2} & = & \s A_{x1}A_{x5}^{-1}\s^{-1} = \s A_{5x}A_{1x}^{-1}\s^{-1}\qquad
     where\ x\ne 1,5 \\
A_{\s,3} & = & \s A_{x2}A_{x3}^{-1}\s^{-1} = \s A_{3x}A_{2x}^{-1}\s^{-1}\qquad
     where\ x\ne 2,3 \\
A_{\s,4} & = & \s A_{x6}A_{x2}^{-1}\s^{-1} = \s A_{2x}A_{6x}^{-1}\s^{-1}\qquad
     where\ x\ne 2,6 \\
A_{\s,5} & = & \s A_{x4}A_{x6}^{-1}\s^{-1} = \s A_{6x}A_{4x}^{-1}\s^{-1}\qquad
     where\ x\ne 4,6 \\
A_{\s,6} & = & \s A_{x5}A_{x4}^{-1}\s^{-1} = \s A_{4x}A_{5x}^{-1}\s^{-1}\qquad
     where\ x\ne 4,5
\end{eqnarray}
\end{tab}
The indices $1$ and $5$ which appear in the formula for $A_{\s,2}$ arise
because $\psi(\G_2)=(15)$.  The conjugations by $\s$ change the indices
as in equation (\ref{As6act}).  Similar for $A_{\s,3}\cdots A_{\s,6}$.

We have reduced the generating set for $\mathcal{A}$ to $\{X_{k
\ell}, A_{k \ell} \}_{k \neq \ell}$ and we know that the subgroup
$K = \sg{X_{k\ell}} \cong \Z^5$.  Now we use the
Reidmeister-Schreier rewriting process to translate all of the
relations of $\Pitil_1$.
From now on, we denote $\bar g = \varphi\psi(g)$ for every $g
\in \Pitil_1$. Using the notation of Subsection \ref{RMS}, $\g(\s,
\G_j) = I$ and $\g(\s,\G_{j'})=A_{\s,j}^{-1}$
for $j \neq 1$.  For $j=1$, $\g(\s, \G_1) = X_\s^{-1}$, and
$\g(\s,\G_{1'}) = B_\s^{-1}$. We begin by translating some of the
relations which involve $\G_1$ but not $\G_{1'}$. These will yield
the relations among the $X_{k \ell}$ which we already know, but the
exercise is useful nonetheless because other elements will satisfy
identical sets of relations. 

$\G_1\G_1\stackrel{\tau}{\longmapsto}
\g(I,\G_1)\g(\overline{\G_1},\G_1) = X_I^{-1}X_{(13)}^{-1} =
X^{-1}_{13} X^{-1}_{31}$, so we deduce that $X_{31} = X^{-1}_{13}$
and conjugating we get
\begin{equation}\label{Xrel1}
X_{\ell k} = X^{-1}_{k \ell}.
\end{equation}

The relations $[\G_1, \G_4]$, $[\G_1, \G_5]$, and $[\G_1, \G_6]$
in $\Pitil_1$ produce the same relations on $\sg{X_{k\ell}}$.

Now we translate the triple relations (\ref{tripeq}) for $i, j$
adjacent. We start for example with $\G_1\G_2\G_1\G_2\G_1\G_2$.
Later we will 
continue with $\G_{1'}\G_2\G_{1'}\G_2\G_{1'}\G_2$ and finish with
$(\G_{1'}\G_1\G_{1'})\G_2(\G_{1'}\G_1\G_{1'})\G_2(\G_{1'}\G_1\G_{1'})\G_2$.
For these two indices, there is no need to use any more
relations from $\G_{(1)}\G_{(2)}\G_{(1)}   =
\G_{(2)}\G_{(1)}\G_{(2)}$, since the Invariance Theorem \ref{th2}
showed that all of these relations were consequences of the three above.

The relation $\G_1\G_2\G_1\G_2\G_1\G_2$ translates through $\tau$
to the expression $$\g(I,\G_1) \g(\overline{\G}_1,\G_2)
\g(\overline{\G_1\G_2},\G_1) \g(\overline{\G_1\G_2\G_1},\G_2)
\g(\overline{\G_1\G_2\G_1\G_2}, \G_1)
\g(\overline{\G_1\G_2\G_1\G_2\G_1},\G_2).$$  But since $\g(\s,
\G_2) = I$ we get $\g(I,\G_1) \g(\overline{\G_1\G_2},\G_1)
\g(\overline{\G_1\G_2\G_1\G_2},\G_1)$, and since
$\overline{\G_1\G_2} = (13)(15) = (135)$ we can further simplify to
$$X_I^{-1}X_{(135)}^{-1}X_{(531)}^{-1} = X^{-1}_{13} X^{-1}_{51}
X^{-1}_{35}.$$ Thus $X_{35}X_{51}X_{13} = 1$ and including all
conjugates
\begin{equation}\label{Xrel2}
X_{k\ell}X_{\ell m}X_{mk} = 1.
\end{equation}

 Similarly the
relation $\G_1\G_3\G_1\G_3\G_1\G_3$ translates through $\tau$ to
the expression 
$$\g(I,\G_1) \g(\overline{\G}_1,\G_3)
\g(\overline{\G_1\G_3},\G_1) \g(\overline{\G_1\G_3\G_1},\G_3)
\g(\overline{\G_1\G_3\G_1\G_3},\G_1)
\g(\overline{\G_1\G_3\G_1\G_3\G_1},\G_3)$$ which equals
$X_I^{-1}X_{(123)}^{-1}X_{(321)}^{-1} = X^{-1}_{13} X^{-1}_{32}
X^{-1}_{21}$. Thus $X_{21}X_{32}X_{13} = 1$, and conjugating we
obtain
\begin{equation}\label{Xrel3}
X_{\ell m}X_{k \ell}X_{mk} = 1.
\end{equation}

  Together the relations \eqs{Xrel1}{Xrel3} show that $\sg{X_{\ell m}}$ is
  generated by the five elements $X_{12},\dots,X_{16}$ which will commute,
so that $\sg{X_{k\ell}} \cong \Z^5$.
These are precisely the relations we expected among the $X_{k\ell}$ and no
more.

  We continue with some of the relations of $\Pitil_1$ which
involve $\G_{1'}$ but not $\G_1$.  These yield identical relations among
the $B_{k\ell}$. 

$\G_{1'}\G_{1'} \stackrel{\tau}{\longmapsto}
\g(I,\G_{1'}) \g(\overline{\G_{1'}},\G_{1'}) =
B_I^{-1}B_{(13)}^{-1} = B_{13}^{-1}B_{31}^{-1}$.  Thus
$B_{31}=B_{13}^{-1}$ and by all conjugations $B_{\ell
k}=B_{k\ell}^{-1}$. 

The relations $[\G_{1'},\G_4]$,
$[\G_{1'},\G_5]$, and $[\G_{1'},\G_6]$ produce the same relations
on $\sg{B_{\s}}$.

The relation
$\G_{1'}\G_2\G_{1'}\G_2\G_{1'}\G_2$ translates through $\tau$ to the
expression
{\small
$$\g(I,\G_{1'}) \g(\overline{\G}_{1'},\G_2) \g(\overline{\G_{1'}\G_2},\G_{1'})
\g(\overline{\G_{1'}\G_2\G_{1'}},\G_2)
\g(\overline{\G_{1'}\G_2\G_{1'}\G_2},\G_{1'})
\g(\overline{\G_{1'}\G_2\G_{1'}\G_2\G_{1'}},\G_2).$$}
But since
$\g(\s, \G_2) = I$ we get $\g(I,\G_{1'})
\g(\overline{\G_{1'}\G_2},\G_{1'})
\g(\overline{\G_{1'}\G_2\G_{1'}\G_2},\G_{1'})$, and since
$\overline{\G_{1'}\G_2} = (13)(15) = (135)$ we further simplify to
$B_I^{-1}B_{(135)}^{-1}B_{(531)}^{-1} =
B_{13}^{-1}B_{51}^{-1}B_{35}^{-1}$. Thus $B_{35}B_{51}B_{13} = 1$
and including all conjugations we have $B_{k\ell}B_{\ell m}B_{mk} = 1$.
Similarly the relation $\G_{1'}\G_3\G_{1'}\G_3\G_{1'}\G_3$
translates through $\tau$ to the expression
$B_{13}^{-1}B_{32}^{-1}B_{21}^{-1}$. Thus $B_{21}B_{32}B_{13} = 1$
and including all conjugations $B_{\ell m}B_{k \ell}B_{mk} = 1$.
By the arguments applied above for the $X_{k\ell}$, we also have
that  $\la B_{k \ell}\ra \cong \Z^5$ generated by $B_{1k}$.

    We finish with the last necessary triple relations.
Note that $\overline{\G_{1'}\G_1\G_{1'}} =
\overline{\G}_1$ and $\tau(\G_{1'}\G_1\G_{1'}) =
B_I^{-1}X_{(13)}^{-1}B_I^{-1} = B_{13}^{-1}X_{31}^{-1}B_{13}^{-1}
= B_{13}^{-1}X_{13}B_{13}^{-1}$. So if we define $C_{k \ell}$ to
be $B_{k\ell}X^{-1}_{k\ell} B_{k\ell} = X_{k \ell} A^2_{k \ell}$
then the additional relations are $C_{\ell k } = C^{-1}_{k \ell}$,
$C_{k\ell}C_{\ell m}C_{mk} = 1$, and $C_{\ell m}C_{k \ell}C_{mk} =
1$. By the arguments above, the $\{C_{k \ell}\}$ generate another
copy $\Z^5 \subset \mathcal{A}$.  In fact for each exponent $n$
the elements $X_{k \ell}A^n_{k \ell} = B_{k \ell} X^{-1}_{k \ell}
\cdots X^{-1}_{k \ell} B_{k \ell}$ or $A^n_{\ell k}= X_{k \ell}
B^{-1}_{k \ell} \cdots B^{-1}_{k \ell} X_{k \ell}$ generate a
subgroup  isomorphic to $\Z^5$.

     The relations computed thus far turn out to be all of the relations
in $\pi_1(X_{\Gal}^{\Aff})$.  Once we show that the remaining
relations translated from $\Pitil_1$ are consequences of the relations
above we will have proven the following theorem:

\begin{theorem}\label{th100}
The fundamental group $\pi_1(X_{\Gal}^{\Aff})$ is generated by
elements $\{X_{ij}, A_{ij}\}$ with the relations
\begin{eqnarray}
X_{ji}A_{ji}^n & =& (X_{ij}A_{ij}^n)^{-1},\\
(X_{ij}A_{ij}^n)(X_{jk}A_{jk}^n)(X_{ki}A_{ki}^n)& =& 1,\\
(X_{jk}A_{jk}^n)(X_{ij}A_{ij}^n)(X_{ki}A_{ki}^n)& =& 1,\\
A_{ij}A_{ik}^{-1} & = & A_{k\ell}A_{j\ell}^{-1}
\end{eqnarray}
for every $n \in \Z$.
\end{theorem}

Before we show that the remaining relations are redundant
we prove that the relations above imply that some of the $A_{k\ell}$
commute.  We shall frequently use the fact that
$X_{k\ell}X_{\ell m}=X_{\ell m}X_{k\ell}=X_{km}$ which is a consequence
of \eqs{Xrel1}{Xrel3}.  $B_{k\ell}$ and $C_{k\ell}$ satisfy this as
well.

\begin{lemma}\label{Acomlem}
In $\pi_1(X_{\Gal}^{\Aff})$ we have $[A_{ij},A_{ik}]=1$ and
$[A_{ji},A_{ki}]=1$ for distinct $i,j,k$.
\end{lemma}

\begin{proof}
Starting with $C_{ki}C_{jk}C_{ij}=1$ and use the definition of $C_{ij}$
to rewrite it as
$(B_{ki}X_{ik}B_{ki})(B_{jk}X_{kj}B_{jk})(B_{ij}X_{ji}B_{ij})=
B_{ki}X_{ik}B_{ji}X_{kj}B_{ik}X_{ji}B_{ij}=$\nl
$(B_{ki}X_{ik})(B_{ji}X_{ij})(X_{ki}B_{ik})(X_{ji}B_{ij})=
A_{ik}^{-1}A_{ij}^{-1}A_{ik}A_{ij}=1$.  Thus the commutator
$[A_{ij},A_{ik}]=1$.  The relation (\ref{Akl:rel}) can be used to show
that $[A_{ji},A_{ki}]=1$ as well.
\end{proof}

Now we treat the remainder of the relations in $\Pitil_1$.  For $j\ne 1$
the relation $\G_j\G_j$ translates immediately to the null relation.
Next consider $\G_{j'}\G_{j'}$.

$\G_{3'}\G_{3'}\stackrel{\tau}{\longmapsto}A_{I,3}^{-1}A_{(23),3}^{-1}$.
But taking the inverse and using Table \ref{Vtrans} we get
$A_{(23),3}A_{I,3}=
(A_{x3}A_{x2}^{-1})(A_{x2}A_{x3}^{-1})$ which cancels completely.
Identical computations treat all other values of $j$.

$\G_1\G_{4'}\G_1\G_{4'}\stackrel{\tau}{\longmapsto}
X_I^{-1}A_{(13),4}^{-1}X_{(13)(26)}^{-1}A_{(26),4}^{-1}$ and taking the
inverse we get\nl
$A_{(26),4}X_{(13)(26)}A_{(13),4}X_I$.  Using Table \ref{Vtrans} we have
$A_{12}A_{16}^{-1}X_{31}A_{36}A_{32}^{-1}X_{13}=
(X_{21}B_{12})(B_{61}X_{16})X_{31}(X_{63}B_{36})(B_{23}X_{32})X_{13}=
X_{21}B_{62}B_{26}X_{12}=1$.  Similar calculations show that
$[\G_1,\G_{5'}]=1$ and $[\G_1,\G_{6'}]=1$ are also redundant.

$\G_{1'}\G_{4'}\G_{1'}\G_{4'}\stackrel{\tau}{\longmapsto}
B_I^{-1}A_{(13),4}^{-1}B_{(13)(26)}^{-1}A_{(26),4}^{-1}$ and taking the
inverse we get\nl
$A_{(26),4}B_{(13)(26)}A_{(13),4}B_I$.  Using Table \ref{Vtrans} we have
$A_{12}A_{16}^{-1}B_{31}A_{36}A_{32}^{-1}B_{13}$.\break  Now, by
Lemma \ref{Acomlem} we can commute $A_{12}$ and $A_{16}$ as well as $A_{36}$
and $A_{32}$ to get
$A_{16}^{-1}A_{12}B_{31}A_{32}^{-1}A_{36}B_{13}=
(B_{61}X_{16})(X_{21}B_{12})B_{31}(B_{23}X_{32})(X_{63}B_{36})B_{13}=
B_{61}X_{26}X_{62}B_{16}=1$.  Again, similar calculations work for
$[\G_{1'},\G_{5'}]=1$ and $[\G_{1'},\G_{6'}]=1$.

For non-adjacent $i,j\ne 1$ the relation $[\G_i,\G_j]$ translates directly
to the null relation.  So next we treat $[\G_{i'},\G_j]$.

$\G_{2'}\G_3\G_{2'}\G_3\stackrel{\tau}{\longmapsto}
A_{I,2}^{-1}A_{(15)(23),2}^{-1}$ and inverting we have
$A_{(15)(23),2}A_{I,3}=$\nl
$(A_{x5}A_{x1}^{-1})(A_{x1}A_{x5}^{-1})=1$.  The same
happens for every other non-adjacent pair $i,j\ne 1$.

$\G_{2'}\G_{3'}\G_{2'}\G_{3'}\stackrel{\tau}{\longmapsto}
A_{I,2}^{-1}A_{(15),3}^{-1}A_{(15)(23),2}^{-1}A_{(23),3}^{-1}$ and
taking inverses again we get
$A_{(23),3}A_{(15)(23),2}A_{(15),3}A_{I,2}=
(A_{x3}A_{x2}^{-1})(A_{1y}A_{5y}^{-1})(A_{z2}A_{z3}^{-1})(A_{5w}A_{1w}^{-1})$.
Substituting specific values $x=1$, $y=2$, $z=5$, and $w=3$ we get\nl
$A_{13}A_{12}^{-1}A_{12}A_{52}^{-1}A_{52}A_{53}^{-1}A_{53}A_{13}^{-1}=1$.
Identical arguments work for every other non-adjacent $i,j\ne 1$.

All that remains are the triple relations for $i,j\ne 1$.  As before we
need only three such relations for each pair of indices.  The relation
$\G_i\G_j\G_i\G_j\G_i\G_j$ translates trivially, so we begin with
$\G_{i'}\G_j\G_{i'}\G_j\G_{i'}\G_j$.

$\G_{4'}\G_3\G_{4'}\G_3\G_{4'}\G_3\stackrel{\tau}{\longmapsto}
A_{I,4}^{-1}A_{(632),4}^{-1}A_{(236),4}^{-1}$ and taking the inverse we get\nl
$A_{(236),4}A_{(632),4}A_{I,4}=
(A_{6x}A_{3x}^{-1})(A_{3x}A_{2x}^{-1})(A_{2x}A_{6x}^{-1})=1$.

Finally consider
$\G_4\G_{4'}\G_4\G_3\G_4\G_{4'}\G_4\G_3\G_4\G_{4'}\G_4\G_3
\stackrel{\tau}{\longmapsto}
A_{(26),4}^{-1}A_{(36),4}^{-1}A_{(23),4}^{-1}$ and
$A_{(23),4}A_{(36),4}A_{(26),4}=
(A_{x6}A_{x3}^{-1})(A_{x3}A_{x2}^{-1})(A_{x2}A_{x6}^{-1})=1$.
So all of the relations in $\Pitil_1$ are included in Theorem \ref{th100}.

\section{The Projective Relation}\label{sec:9}

     To complete the computation of $\pi_1(X_{\Gal})$ we need only to add
the projective relation
$$\G_1\G_{1'}\G_2\G_{2'}\G_3\G_{3'}\G_4\G_{4'}\G_5\G_{5'}\G_6\G_{6'} = 1.$$
This relation translates in $\mathcal{A}$ as the product $P =
A_{I,1}A_{I,2}A_{I,3}A_{I,4}A_{I,5}A_{I,6}$.
We must translate the $A_{I,j}$ to the language of the $A_{k
\ell}$, using Table \ref{Vtrans}:

$P$ translates to
$A_{13}(A_{21}A_{25}^{-1})(A_{31}A_{21}^{-1})(A_{21}A_{61}^{-1})
(A_{61}A_{41}^{-1})(A_{41}A_{51}^{-1})$ which cancels to
$A_{13}A_{21}A_{25}^{-1}A_{31}A_{51}^{-1}$.  Using Equation
\eq{Akl:rel}, we get $A_{13}A_{21}A_{25}^{-1}A_{25}A_{23}^{-1}$.
Thus the projective relation may be written as
$A_{13}A_{21}A_{23}^{-1}=1$ or equivalently $A_{23}=A_{13}A_{21}$.
Conjugating, this becomes
\begin{equation}\label{Akl:rel2}
A_{ij} = A_{kj}A_{ik}.
\end{equation}
Substituting back into \eq{Akl:rel}, writing $A_{ij} =
A_{kj}A_{ik}$ and $A_{k\ell} = A_{j\ell}A_{kj}$, we obtain
\begin{equation}\label{Akl:rel3}
A_{kj} A_{j \ell} = A_{j\ell}A_{kj}.
\end{equation}

\begin{lemma}
The subgroup $\sg{A_{k\ell}}$ of $\pi_1(X_{\Gal})$ is
commutative of rank of at most $5$
\end{lemma}
\begin{proof}
We will compute the centralizer of $A_{ij}$ for fixed $i,j$. Let
$i,j,k,\ell$ be four distinct indices.  We already know from Lemma
\ref{Acomlem} that $A_{ij}$ commutes with $A_{ik}$ and $A_{\ell j}$.
By equation (\ref{Akl:rel3}) it also commutes with $A_{ki}$ and
$A_{j\ell}$.  Now equation (\ref{Akl:rel2}) allows us to write
$A_{k\ell}=A_{i\ell}A_{\ell j}$, both of which commute with $A_{ij}$,
so $\sg{A_{k\ell}}$ is commutative.

Now, since $A_{j k} A_{ij} = A_{ik}$ and $A_{ik}A_{ji} = A_{jk}$, we have
$A_{jk} A_{ij} A_{ji} = A_{ik} A_{ji} = A_{jk}$, so that $A_{ji} =
A_{ij}^{-1}$, the group is generated by the $A_{1k}$ ($k = 2,\dots,6$), and
the rank is at most $5$.
\end{proof}

We see that $\pi_1(X_{\Gal}) = \sg{A_{ij},X_{ij}}$ with the two
  subgroups $\sg{A_{ij}},\sg{X_{ij}}$ isomorphic to  $\Z^5$.
   The only question left is how these
two subgroups interact.

\begin{lemma}\label{lm101}
In $\pi_1(X_{\Gal})$ the $A_{ij}$ and $X_{k\ell}$ commute.
\end{lemma}

\begin{proof}
We need only consider the commutators of $A_{13}$ and $X_{ij}$
since all others are merely conjugates of these. First consider
the commutator $[X_{13},A_{13}]$. Since
$X_{13}=(13)\G_1$ (choose $\s = 1$ in \eq{Xkl} and note that as elements of 
$S_6$, we have
$(13) = \G_2\G_6\G_5\G_4\G_3\G_4\G_5\G_6\G_2$)
and $A_{13}=\G_1\G_{1'}$ this becomes
$(13)\G_1(\G_1\G_{1'})\G_1(13)A_{13}^{-1}=(13)\G_{1'}\G_1(13)A_{13}^{-1}=
A_{31}^{-1}A_{13}^{-1}=A_{13}A_{13}^{-1}=1$.  So $X_{13}$ and
$A_{13}$ commute.

Next consider the commutator $X_{12}A_{13}X_{12}^{-1}A_{13}^{-1}$.
By definition we have that 
$X_{12}=(23)X_{13}(23)=(23)(13)\G_1(23)=(321)\G_1\G_3$. Thus the
commutator becomes
$(321)\G_1\G_3(\G_1\G_{1'})\G_3\G_1(123)A_{13}^{-1}$. We use the
triple relations\nl $\Trip{\Gamma_{(1)}}{\Gamma_{(3)}}$ to get
$(321)\G_3\G_1\G_3\G_{1'}\G_3\G_1(123)A_{13}^{-1}=$\nl$
(321)\G_3\G_1\G_{1'}\G_3\G_{1'}\G_1(123)A_{13}^{-1}=(321)\G_3\G_1\G_{1'}\G_3(123)(321)\G_{1'}\G_1(123)A_{13}^{-1}$\nl
which is equal to $A_{(321)\G_3,1}A_{(321),1}^{-1}A_{13}^{-1}=
A_{(12),1}A_{(321),1}^{-1}A_{13}^{-1}=A_{23}A_{21}^{-1}A_{13}^{-1}=1$, 
thus proving that $X_{12}$ and $A_{13}$ commute.

Conjugating by $(2j)$ we see that $X_{1j}$ commutes with $A_{13}$
and since $X_{ij}=X_{1i}^{-1}X_{1j}$ we see that every $X_{ij}$
commutes with $A_{13}$.
\end{proof}

\begin{theorem}\label{main}
The fundamental group $\pi_1(X_{\Gal})\cong\Z^{10}$.
\end{theorem}
\begin{proof}
$\pi_1(X_{\Gal})$ is generated by $A_{1j}$ and $X_{1j}$ which all
commute.  Hence the group they generate is $\Z^{10}$.
\end{proof}

\Addresses\recd

\end{document}

%% file: agtout.tex
\def\ifplaintex{\expandafter\ifx\csname documentclass\endcsname\relax}

\def\gtp{{\mathsurround=0pt\it $\cal G\mskip-2mu$eometry \&\ 
$\cal T\!\!$opology $\cal P\!$ublications}}  

\def\recd{{\small Received:\qua\receiveddate\ifx\reviseddate\relax
\else\qquad Revised:\qua\reviseddate\fi\par}} 

\def\lognumber#1{\def\thelognumber{#1}}
\def\volumenumber#1{\def\thevolumenumber{#1}}
\def\volumeyear#1{\def\thevolumeyear{#1}}
\def\papernumber#1{\def\thepapernumber{#1}}
\def\pagenumbers#1#2{\def\startpage{#1}\def\finishpage{#2}}
\def\published#1{\def\publishdate{#1}}

\def\received#1{\def\receiveddate{#1}}
\def\revised#1{\def\reviseddate{#1}}
\def\accepted#1{\def\accepteddate{#1}}
\def\asciititle#1{\def\theasciititle{#1}}
\def\covertitle#1{\def\thecovertitle{#1}}

\def\asciiaddress#1{\def\theasciiaddress{#1}}

\long\def\asciiabstract#1{\long\def\theasciiabstract{#1}}

\let\\\par\let\thelognumber\relax\let\thevolumenumber\relax
\let\thepapernumber\relax\let\thevolumeyear\relax\let\startpage\relax
\let\finishpage\relax\let\publishdate\relax\let\receiveddate\relax
\let\reviseddate\relax\let\accepteddate\relax\let\theasciititle\relax
\let\thecovertitle\relax\let\theasciiauthors\relax\let\theasciiaddress\relax
\let\theasciiabstract\relax

\let\theasciiemail\relax

\ifplaintex
\font\logobig=cmssbx10 scaled 3836
\font\logomed=cmssbx10 scaled 2557
\else
\font\logobig=cmssbx10 scaled 4200
\font\logomed=cmssbx10 scaled 2800
\fi

\long\def\makeagttitle{   
\count0=\startpage
\agt\hfill      
\hbox to 45truept{\vbox to 0pt{\vglue -13truept{\logomed A\kern -.37em{\logobig 
T}\kern -.38em G}\vss}\hss}
\break
{\small Volume \thevolumenumber\ (\thevolumeyear)
\startpage--\finishpage\nl
Published: \publishdate}

\vglue .25truein

{\parskip=0pt\leftskip 0pt plus
1fil\def\\{\par\smallskip}{\Large\bf\thetitle}\par\medskip} \vglue
0.05truein

{\parskip=0pt\leftskip 0pt plus 1fil\def\\{\par}{\sc\theauthors}
\par\medskip}%
 
\vglue 0.03truein 

{\small\leftskip 25truept\rightskip 25truept{\bf Abstract}\stdspace\theabstract

{\bf AMS Classification}\stdspace\theprimaryclass
\ifx\thesecondaryclass\relax\else; \thesecondaryclass\fi\par
{\bf Keywords}\stdspace \thekeywords\par}\vglue 7truept

}   

\ifplaintex
\font\phead=cmsl9 scaled 950
\font\pnum=cmbx10 scaled 913
\font\pfoot=cmsl9 scaled 950
\headline{\vbox to 0pt{\vskip -4.5mm\line{\small\phead\ifnum
\count0=\startpage ISSN 1472-2739 (on-line) 1472-2747 (printed)
\hfill {\pnum\folio}\else\ifodd\count0\def\\{ }%
\ifx\theshorttitle\relax\thetitle\else\theshorttitle\fi\hfill{\pnum\folio}
\else\def\\{ and }{\pnum\folio}\hfill\ifx\theshortauthors\relax\theauthors
\else\theshortauthors\fi\fi\fi}\vss}}
\footline{\vbox to 0pt{\vglue 0mm\line{\small\pfoot\ifnum\count0=\startpage
\copyright\ \gtp\hfill\else
\agt, Volume \thevolumenumber\ (\thevolumeyear)\hfill\fi}\vss}}
\else
\font=cmsl9 scaled 1050
\font\lnum=cmbx10 
\font=cmsl9 scaled 1050
\makeatletter
\def\@oddhead{{\small\lhead\ifnum\count0=\startpage ISSN 1472-2739 
(on-line) 1472-2747 (printed)\hfill {\lnum\number\count0}\else\ifodd\count0
\def\\{ }\ifx\theshorttitle\relax \thetitle \else\theshorttitle\fi\hfill
{\lnum\number\count0}\else\def\\{ and }{\lnum\number\count0}
\hfill\ifx\theshortauthors\relax 
\theauthors\else\theshortauthors\fi\fi\fi}}\def\@evenhead{\@oddhead}
\def\@oddfoot{\small\lfoot\ifnum\count0=\startpage\copyright\ \gtp\hfill\else
\agt, Volume \thevolumenumber\ (\thevolumeyear)\hfill\fi}
\def\@evenfoot{\@oddfoot}
\makeatother
\fi
\let\maketitlepage\makeagttitle
\let\makeshorttitle\maketitlepage
\let\maketitle\maketitlepage

\newwrite\gtoutfile
\long\gdef\makeheadfile{  
{\def\\{, }\def\s{ }
\immediate\openout\gtoutfile head.xxx
\immediate\write\gtoutfile{To: math@arxiv.org}
\immediate\write\gtoutfile{Subject: put OR rep NNNNN:ppppp}
\immediate\write\gtoutfile{--text follows this line--}
\immediate\write\gtoutfile{Proxy-for: \ifx\theasciiauthors\relax
\theauthors\else\theasciiauthors\fi\s<\ifx\theasciiemail\relax\theemail\else\theasciiemail\fi>}
\immediate\write\gtoutfile{\noexpand\\}
\immediate\write\gtoutfile{Authors: \ifx\theasciiauthors\relax
\theauthors\else\theasciiauthors\fi}
{\def\\{ }\immediate\write\gtoutfile{Title: \ifx\theasciititle\relax
\thetitle\else\theasciititle\fi}}
\immediate\write\gtoutfile{Subj-class: GT or SG, GR etc}
\immediate\write\gtoutfile{MSC-class: \theprimaryclass\ifx\thesecondaryclass\relax\else, \thesecondaryclass\fi}
\immediate\write\gtoutfile{Journal-ref: Algebr. Geom. Topol. \thevolumenumber\s
(\thevolumeyear) \startpage-\finishpage}
\immediate\write\gtoutfile{Comments: Published by Algebraic and
Geometric Topology at}
\immediate\write\gtoutfile{\s\s\s  http://www.maths.warwick.ac.uk/agt/AGTVol\thevolumenumber/agt-\thevolumenumber-\thepapernumber.abs.html}
\immediate\write\gtoutfile{\noexpand\\}
\immediate\write\gtoutfile{}
\ifx\theasciiabstract\relax
\immediate\write\gtoutfile{\theabstract}\else
\immediate\write\gtoutfile{\theasciiabstract}\fi
\immediate\write\gtoutfile{}
\immediate\write\gtoutfile{\noexpand\\}
\immediate\write\gtoutfile{}
\immediate\closeout\gtoutfile}}  

\def\maketitlepage{\makeagttitle\makeheadfile}
\let\makeshorttitle\maketitlepage
\let\maketitle\maketitlepage

\def\ifplaintex{\expandafter\ifx\csname documentclass\endcsname\relax}

\def\gtp{{\mathsurround=0pt\it $\cal G\mskip-2mu$eometry \&\ 
$\cal T\!\!$opology $\cal P\!$ublications}}  

\def\recd{{\small Received:\qua\receiveddate\ifx\reviseddate\relax
\else\qquad Revised:\qua\reviseddate\fi\par}} 

\def\lognumber#1{\def\thelognumber{#1}}
\def\volumenumber#1{\def\thevolumenumber{#1}}
\def\volumeyear#1{\def\thevolumeyear{#1}}
\def\papernumber#1{\def\thepapernumber{#1}}
\def\pagenumbers#1#2{\def\startpage{#1}\def\finishpage{#2}}
\def\published#1{\def\publishdate{#1}}

\def\received#1{\def\receiveddate{#1}}
\def\revised#1{\def\reviseddate{#1}}
\def\accepted#1{\def\accepteddate{#1}}
\def\asciititle#1{\def\theasciititle{#1}}
\def\covertitle#1{\def\thecovertitle{#1}}

\def\asciiaddress#1{\def\theasciiaddress{#1}}

\long\def\asciiabstract#1{\long\def\theasciiabstract{#1}}

\let\\\par\let\thelognumber\relax\let\thevolumenumber\relax
\let\thepapernumber\relax\let\thevolumeyear\relax\let\startpage\relax
\let\finishpage\relax\let\publishdate\relax\let\receiveddate\relax
\let\reviseddate\relax\let\accepteddate\relax\let\theasciititle\relax
\let\thecovertitle\relax\let\theasciiauthors\relax\let\theasciiaddress\relax
\let\theasciiabstract\relax

\let\theasciiemail\relax

\ifplaintex
\font\logobig=cmssbx10 scaled 3836
\font\logomed=cmssbx10 scaled 2557
\else
\font\logobig=cmssbx10 scaled 4200
\font\logomed=cmssbx10 scaled 2800
\fi

\long\def\makeagttitle{   
\count0=\startpage
\agt\hfill      
\hbox to 45truept{\vbox to 0pt{\vglue -13truept{\logomed A\kern -.37em{\logobig 
T}\kern -.38em G}\vss}\hss}
\break
{\small Volume \thevolumenumber\ (\thevolumeyear)
\startpage--\finishpage\nl
Published: \publishdate}

\vglue .25truein

{\parskip=0pt\leftskip 0pt plus
1fil\def\\{\par\smallskip}{\Large\bf\thetitle}\par\medskip} \vglue
0.05truein

{\parskip=0pt\leftskip 0pt plus 1fil\def\\{\par}{\sc\theauthors}
\par\medskip}%
 
\vglue 0.03truein 

{\small\leftskip 25truept\rightskip 25truept{\bf Abstract}\stdspace\theabstract

{\bf AMS Classification}\stdspace\theprimaryclass
\ifx\thesecondaryclass\relax\else; \thesecondaryclass\fi\par
{\bf Keywords}\stdspace \thekeywords\par}\vglue 7truept

}   

\ifplaintex
\font\phead=cmsl9 scaled 950
\font\pnum=cmbx10 scaled 913
\font\pfoot=cmsl9 scaled 950
\headline{\vbox to 0pt{\vskip -4.5mm\line{\small\phead\ifnum
\count0=\startpage ISSN 1472-2739 (on-line) 1472-2747 (printed)
\hfill {\pnum\folio}\else\ifodd\count0\def\\{ }%
\ifx\theshorttitle\relax\thetitle\else\theshorttitle\fi\hfill{\pnum\folio}
\else\def\\{ and }{\pnum\folio}\hfill\ifx\theshortauthors\relax\theauthors
\else\theshortauthors\fi\fi\fi}\vss}}
\footline{\vbox to 0pt{\vglue 0mm\line{\small\pfoot\ifnum\count0=\startpage
\copyright\ \gtp\hfill\else
\agt, Volume \thevolumenumber\ (\thevolumeyear)\hfill\fi}\vss}}
\else
\font=cmsl9 scaled 1050
\font\lnum=cmbx10 
\font=cmsl9 scaled 1050
\makeatletter
\def\@oddhead{{\small\lhead\ifnum\count0=\startpage ISSN 1472-2739 
(on-line) 1472-2747 (printed)\hfill {\lnum\number\count0}\else\ifodd\count0
\def\\{ }\ifx\theshorttitle\relax \thetitle \else\theshorttitle\fi\hfill
{\lnum\number\count0}\else\def\\{ and }{\lnum\number\count0}
\hfill\ifx\theshortauthors\relax 
\theauthors\else\theshortauthors\fi\fi\fi}}\def\@evenhead{\@oddhead}
\def\@oddfoot{\small\lfoot\ifnum\count0=\startpage\copyright\ \gtp\hfill\else
\agt, Volume \thevolumenumber\ (\thevolumeyear)\hfill\fi}
\def\@evenfoot{\@oddfoot}
\makeatother
\fi
\let\maketitlepage\makeagttitle
\let\makeshorttitle\maketitlepage
\let\maketitle\maketitlepage

\newwrite\gtoutfile
\long\gdef\makeheadfile{  
{\def\\{, }\def\s{ }
\immediate\openout\gtoutfile head.xxx
\immediate\write\gtoutfile{To: math@arxiv.org}
\immediate\write\gtoutfile{Subject: put OR rep NNNNN:ppppp}
\immediate\write\gtoutfile{--text follows this line--}
\immediate\write\gtoutfile{Proxy-for: \ifx\theasciiauthors\relax
\theauthors\else\theasciiauthors\fi\s<\ifx\theasciiemail\relax\theemail\else\theasciiemail\fi>}
\immediate\write\gtoutfile{\noexpand\\}
\immediate\write\gtoutfile{Authors: \ifx\theasciiauthors\relax
\theauthors\else\theasciiauthors\fi}
{\def\\{ }\immediate\write\gtoutfile{Title: \ifx\theasciititle\relax
\thetitle\else\theasciititle\fi}}
\immediate\write\gtoutfile{Subj-class: GT or SG, GR etc}
\immediate\write\gtoutfile{MSC-class: \theprimaryclass\ifx\thesecondaryclass\relax\else, \thesecondaryclass\fi}
\immediate\write\gtoutfile{Journal-ref: Algebr. Geom. Topol. \thevolumenumber\s
(\thevolumeyear) \startpage-\finishpage}
\immediate\write\gtoutfile{Comments: Published by Algebraic and
Geometric Topology at}
\immediate\write\gtoutfile{\s\s\s  http://www.maths.warwick.ac.uk/agt/AGTVol\thevolumenumber/agt-\thevolumenumber-\thepapernumber.abs.html}
\immediate\write\gtoutfile{\noexpand\\}
\immediate\write\gtoutfile{}
\ifx\theasciiabstract\relax
\immediate\write\gtoutfile{\theabstract}\else
\immediate\write\gtoutfile{\theasciiabstract}\fi
\immediate\write\gtoutfile{}
\immediate\write\gtoutfile{\noexpand\\}
\immediate\write\gtoutfile{}
\immediate\closeout\gtoutfile}}  

\def\maketitlepage{\makeagttitle\makeheadfile}
\let\makeshorttitle\maketitlepage
\let\maketitle\maketitlepage

\def\ifplaintex{\expandafter\ifx\csname documentclass\endcsname\relax}

\def\gtp{{\mathsurround=0pt\it $\cal G\mskip-2mu$eometry \&\ 
$\cal T\!\!$opology $\cal P\!$ublications}}  

\def\recd{{\small Received:\qua\receiveddate\ifx\reviseddate\relax
\else\qquad Revised:\qua\reviseddate\fi\par}} 

\def\lognumber#1{\def\thelognumber{#1}}
\def\volumenumber#1{\def\thevolumenumber{#1}}
\def\volumeyear#1{\def\thevolumeyear{#1}}
\def\papernumber#1{\def\thepapernumber{#1}}
\def\pagenumbers#1#2{\def\startpage{#1}\def\finishpage{#2}}
\def\published#1{\def\publishdate{#1}}

\def\received#1{\def\receiveddate{#1}}
\def\revised#1{\def\reviseddate{#1}}
\def\accepted#1{\def\accepteddate{#1}}
\def\asciititle#1{\def\theasciititle{#1}}
\def\covertitle#1{\def\thecovertitle{#1}}

\def\asciiaddress#1{\def\theasciiaddress{#1}}

\long\def\asciiabstract#1{\long\def\theasciiabstract{#1}}

\let\\\par\let\thelognumber\relax\let\thevolumenumber\relax
\let\thepapernumber\relax\let\thevolumeyear\relax\let\startpage\relax
\let\finishpage\relax\let\publishdate\relax\let\receiveddate\relax
\let\reviseddate\relax\let\accepteddate\relax\let\theasciititle\relax
\let\thecovertitle\relax\let\theasciiauthors\relax\let\theasciiaddress\relax
\let\theasciiabstract\relax

\let\theasciiemail\relax

\ifplaintex
\font\logobig=cmssbx10 scaled 3836
\font\logomed=cmssbx10 scaled 2557
\else
\font\logobig=cmssbx10 scaled 4200
\font\logomed=cmssbx10 scaled 2800
\fi

\long\def\makeagttitle{   
\count0=\startpage
\agt\hfill      
\hbox to 45truept{\vbox to 0pt{\vglue -13truept{\logomed A\kern -.37em{\logobig 
T}\kern -.38em G}\vss}\hss}
\break
{\small Volume \thevolumenumber\ (\thevolumeyear)
\startpage--\finishpage\nl
Published: \publishdate}

\vglue .25truein

{\parskip=0pt\leftskip 0pt plus
1fil\def\\{\par\smallskip}{\Large\bf\thetitle}\par\medskip} \vglue
0.05truein

{\parskip=0pt\leftskip 0pt plus 1fil\def\\{\par}{\sc\theauthors}
\par\medskip}%
 
\vglue 0.03truein 

{\small\leftskip 25truept\rightskip 25truept{\bf Abstract}\stdspace\theabstract

{\bf AMS Classification}\stdspace\theprimaryclass
\ifx\thesecondaryclass\relax\else; \thesecondaryclass\fi\par
{\bf Keywords}\stdspace \thekeywords\par}\vglue 7truept

}   

\ifplaintex
\font\phead=cmsl9 scaled 950
\font\pnum=cmbx10 scaled 913
\font\pfoot=cmsl9 scaled 950
\headline{\vbox to 0pt{\vskip -4.5mm\line{\small\phead\ifnum
\count0=\startpage ISSN 1472-2739 (on-line) 1472-2747 (printed)
\hfill {\pnum\folio}\else\ifodd\count0\def\\{ }%
\ifx\theshorttitle\relax\thetitle\else\theshorttitle\fi\hfill{\pnum\folio}
\else\def\\{ and }{\pnum\folio}\hfill\ifx\theshortauthors\relax\theauthors
\else\theshortauthors\fi\fi\fi}\vss}}
\footline{\vbox to 0pt{\vglue 0mm\line{\small\pfoot\ifnum\count0=\startpage
\copyright\ \gtp\hfill\else
\agt, Volume \thevolumenumber\ (\thevolumeyear)\hfill\fi}\vss}}
\else
\font=cmsl9 scaled 1050
\font\lnum=cmbx10 
\font=cmsl9 scaled 1050
\makeatletter
\def\@oddhead{{\small\lhead\ifnum\count0=\startpage ISSN 1472-2739 
(on-line) 1472-2747 (printed)\hfill {\lnum\number\count0}\else\ifodd\count0
\def\\{ }\ifx\theshorttitle\relax \thetitle \else\theshorttitle\fi\hfill
{\lnum\number\count0}\else\def\\{ and }{\lnum\number\count0}
\hfill\ifx\theshortauthors\relax 
\theauthors\else\theshortauthors\fi\fi\fi}}\def\@evenhead{\@oddhead}
\def\@oddfoot{\small\lfoot\ifnum\count0=\startpage\copyright\ \gtp\hfill\else
\agt, Volume \thevolumenumber\ (\thevolumeyear)\hfill\fi}
\def\@evenfoot{\@oddfoot}
\makeatother
\fi
\let\maketitlepage\makeagttitle
\let\makeshorttitle\maketitlepage
\let\maketitle\maketitlepage

\newwrite\gtoutfile
\long\gdef\makeheadfile{  
{\def\\{, }\def\s{ }
\immediate\openout\gtoutfile head.xxx
\immediate\write\gtoutfile{To: math@arxiv.org}
\immediate\write\gtoutfile{Subject: put OR rep NNNNN:ppppp}
\immediate\write\gtoutfile{--text follows this line--}
\immediate\write\gtoutfile{Proxy-for: \ifx\theasciiauthors\relax
\theauthors\else\theasciiauthors\fi\s<\ifx\theasciiemail\relax\theemail\else\theasciiemail\fi>}
\immediate\write\gtoutfile{\noexpand\\}
\immediate\write\gtoutfile{Authors: \ifx\theasciiauthors\relax
\theauthors\else\theasciiauthors\fi}
{\def\\{ }\immediate\write\gtoutfile{Title: \ifx\theasciititle\relax
\thetitle\else\theasciititle\fi}}
\immediate\write\gtoutfile{Subj-class: GT or SG, GR etc}
\immediate\write\gtoutfile{MSC-class: \theprimaryclass\ifx\thesecondaryclass\relax\else, \thesecondaryclass\fi}
\immediate\write\gtoutfile{Journal-ref: Algebr. Geom. Topol. \thevolumenumber\s
(\thevolumeyear) \startpage-\finishpage}
\immediate\write\gtoutfile{Comments: Published by Algebraic and
Geometric Topology at}
\immediate\write\gtoutfile{\s\s\s  http://www.maths.warwick.ac.uk/agt/AGTVol\thevolumenumber/agt-\thevolumenumber-\thepapernumber.abs.html}
\immediate\write\gtoutfile{\noexpand\\}
\immediate\write\gtoutfile{}
\ifx\theasciiabstract\relax
\immediate\write\gtoutfile{\theabstract}\else
\immediate\write\gtoutfile{\theasciiabstract}\fi
\immediate\write\gtoutfile{}
\immediate\write\gtoutfile{\noexpand\\}
\immediate\write\gtoutfile{}
\immediate\closeout\gtoutfile}}  

\def\maketitlepage{\makeagttitle\makeheadfile}
\let\makeshorttitle\maketitlepage
\let\maketitle\maketitlepage

\def\ifplaintex{\expandafter\ifx\csname documentclass\endcsname\relax}

\def\gtp{{\mathsurround=0pt\it $\cal G\mskip-2mu$eometry \&\ 
$\cal T\!\!$opology $\cal P\!$ublications}}  

\def\recd{{\small Received:\qua\receiveddate\ifx\reviseddate\relax
\else\qquad Revised:\qua\reviseddate\fi\par}} 

\def\lognumber#1{\def\thelognumber{#1}}
\def\volumenumber#1{\def\thevolumenumber{#1}}
\def\volumeyear#1{\def\thevolumeyear{#1}}
\def\papernumber#1{\def\thepapernumber{#1}}
\def\pagenumbers#1#2{\def\startpage{#1}\def\finishpage{#2}}
\def\published#1{\def\publishdate{#1}}

\def\received#1{\def\receiveddate{#1}}
\def\revised#1{\def\reviseddate{#1}}
\def\accepted#1{\def\accepteddate{#1}}
\def\asciititle#1{\def\theasciititle{#1}}
\def\covertitle#1{\def\thecovertitle{#1}}

\def\asciiaddress#1{\def\theasciiaddress{#1}}

\long\def\asciiabstract#1{\long\def\theasciiabstract{#1}}

\let\\\par\let\thelognumber\relax\let\thevolumenumber\relax
\let\thepapernumber\relax\let\thevolumeyear\relax\let\startpage\relax
\let\finishpage\relax\let\publishdate\relax\let\receiveddate\relax
\let\reviseddate\relax\let\accepteddate\relax\let\theasciititle\relax
\let\thecovertitle\relax\let\theasciiauthors\relax\let\theasciiaddress\relax
\let\theasciiabstract\relax

\let\theasciiemail\relax

\ifplaintex
\font\logobig=cmssbx10 scaled 3836
\font\logomed=cmssbx10 scaled 2557
\else
\font\logobig=cmssbx10 scaled 4200
\font\logomed=cmssbx10 scaled 2800
\fi

\long\def\makeagttitle{   
\count0=\startpage
\agt\hfill      
\hbox to 45truept{\vbox to 0pt{\vglue -13truept{\logomed A\kern -.37em{\logobig 
T}\kern -.38em G}\vss}\hss}
\break
{\small Volume \thevolumenumber\ (\thevolumeyear)
\startpage--\finishpage\nl
Published: \publishdate}

\vglue .25truein

{\parskip=0pt\leftskip 0pt plus
1fil\def\\{\par\smallskip}{\Large\bf\thetitle}\par\medskip} \vglue
0.05truein

{\parskip=0pt\leftskip 0pt plus 1fil\def\\{\par}{\sc\theauthors}
\par\medskip}%
 
\vglue 0.03truein 

{\small\leftskip 25truept\rightskip 25truept{\bf Abstract}\stdspace\theabstract

{\bf AMS Classification}\stdspace\theprimaryclass
\ifx\thesecondaryclass\relax\else; \thesecondaryclass\fi\par
{\bf Keywords}\stdspace \thekeywords\par}\vglue 7truept

}   

\ifplaintex
\font\phead=cmsl9 scaled 950
\font\pnum=cmbx10 scaled 913
\font\pfoot=cmsl9 scaled 950
\headline{\vbox to 0pt{\vskip -4.5mm\line{\small\phead\ifnum
\count0=\startpage ISSN 1472-2739 (on-line) 1472-2747 (printed)
\hfill {\pnum\folio}\else\ifodd\count0\def\\{ }%
\ifx\theshorttitle\relax\thetitle\else\theshorttitle\fi\hfill{\pnum\folio}
\else\def\\{ and }{\pnum\folio}\hfill\ifx\theshortauthors\relax\theauthors
\else\theshortauthors\fi\fi\fi}\vss}}
\footline{\vbox to 0pt{\vglue 0mm\line{\small\pfoot\ifnum\count0=\startpage
\copyright\ \gtp\hfill\else
\agt, Volume \thevolumenumber\ (\thevolumeyear)\hfill\fi}\vss}}
\else
\font=cmsl9 scaled 1050
\font\lnum=cmbx10 
\font=cmsl9 scaled 1050
\makeatletter
\def\@oddhead{{\small\lhead\ifnum\count0=\startpage ISSN 1472-2739 
(on-line) 1472-2747 (printed)\hfill {\lnum\number\count0}\else\ifodd\count0
\def\\{ }\ifx\theshorttitle\relax \thetitle \else\theshorttitle\fi\hfill
{\lnum\number\count0}\else\def\\{ and }{\lnum\number\count0}
\hfill\ifx\theshortauthors\relax 
\theauthors\else\theshortauthors\fi\fi\fi}}\def\@evenhead{\@oddhead}
\def\@oddfoot{\small\lfoot\ifnum\count0=\startpage\copyright\ \gtp\hfill\else
\agt, Volume \thevolumenumber\ (\thevolumeyear)\hfill\fi}
\def\@evenfoot{\@oddfoot}
\makeatother
\fi
\let\maketitlepage\makeagttitle
\let\makeshorttitle\maketitlepage
\let\maketitle\maketitlepage

\newwrite\gtoutfile
\long\gdef\makeheadfile{  
{\def\\{, }\def\s{ }
\immediate\openout\gtoutfile head.xxx
\immediate\write\gtoutfile{To: math@arxiv.org}
\immediate\write\gtoutfile{Subject: put OR rep NNNNN:ppppp}
\immediate\write\gtoutfile{--text follows this line--}
\immediate\write\gtoutfile{Proxy-for: \ifx\theasciiauthors\relax
\theauthors\else\theasciiauthors\fi\s<\ifx\theasciiemail\relax\theemail\else\theasciiemail\fi>}
\immediate\write\gtoutfile{\noexpand\\}
\immediate\write\gtoutfile{Authors: \ifx\theasciiauthors\relax
\theauthors\else\theasciiauthors\fi}
{\def\\{ }\immediate\write\gtoutfile{Title: \ifx\theasciititle\relax
\thetitle\else\theasciititle\fi}}
\immediate\write\gtoutfile{Subj-class: GT or SG, GR etc}
\immediate\write\gtoutfile{MSC-class: \theprimaryclass\ifx\thesecondaryclass\relax\else, \thesecondaryclass\fi}
\immediate\write\gtoutfile{Journal-ref: Algebr. Geom. Topol. \thevolumenumber\s
(\thevolumeyear) \startpage-\finishpage}
\immediate\write\gtoutfile{Comments: Published by Algebraic and
Geometric Topology at}
\immediate\write\gtoutfile{\s\s\s  http://www.maths.warwick.ac.uk/agt/AGTVol\thevolumenumber/agt-\thevolumenumber-\thepapernumber.abs.html}
\immediate\write\gtoutfile{\noexpand\\}
\immediate\write\gtoutfile{}
\ifx\theasciiabstract\relax
\immediate\write\gtoutfile{\theabstract}\else
\immediate\write\gtoutfile{\theasciiabstract}\fi
\immediate\write\gtoutfile{}
\immediate\write\gtoutfile{\noexpand\\}
\immediate\write\gtoutfile{}
\immediate\closeout\gtoutfile}}  

\def\maketitlepage{\makeagttitle\makeheadfile}
\let\makeshorttitle\maketitlepage
\let\maketitle\maketitlepage